\documentclass[letterpaper]{article}

\usepackage[backend=biber,sorting=nyvt,style=numeric-comp,isbn=false,doi=false,date=year,autolang=none,abbreviate=true]{biblatex}
\AtEveryBibitem{%
\clearfield{note}%
  \clearfield{number}}

\DeclareSourcemap{
  \maps[datatype=bibtex]{
    \map[overwrite]{ 
      \step[fieldsource=shortjournal,fieldtarget=journaltitle]
    }
  }
}

\DeclareFieldFormat[article,periodical]{volume}{\mkbibbold{#1}}

\renewbibmacro{in:}{%
  \ifentrytype{article}{}{\printtext{\bibstring{in}\intitlepunct}}}

\DeclareFieldFormat[article]{pages}{#1}

\DeclareFieldFormat
  [article]
  {title}{#1\isdot}

\DeclareDelimFormat[bib]{nametitledelim}{\addcomma\space}

\DeclareFieldFormat{url}{\url{#1}}

\DeclareFieldFormat[article,online,book,misc]{title}{#1\addcomma}

\usepackage{amsmath}
\usepackage[colorlinks=true]{hyperref}
\usepackage{amsmath}
\usepackage{amssymb}
\usepackage{amsthm}
\usepackage{amsfonts}
\usepackage{graphicx}
\usepackage[left=2.1cm,right=2.1cm,top=3cm,bottom=2.5cm]{geometry}
\usepackage{subcaption}

\newcommand{\revised}[1]{\textcolor{black}{#1}}

\hypersetup{urlcolor=blue, citecolor=red}

\title{Analysis of Iterative Deblurring: No Explicit Noise } 

\author{%
	Sinethemba Neliswa Mamba\textsuperscript{1} and Pawe\l\ Danielewicz\textsuperscript{2,1}\thanks{Corresponding author: \href{mailto:daniel@frib.msu.edu}{daniel@frib.msu.edu}} \footnote{Pawe\l\ Danielewicz was partly supported by the U.S.\ Department of Energy Office of Science Grant No.\ DE-SC0019209.}
}

\date{\footnotesize\textsuperscript{\textbf{1}}African Institute for Mathematical Sciences, Kigali, Rwanda\\[1ex] \textsuperscript{\textbf{2}}Facility for Rare Isotope Beams and Department of Physics and Astronomy,\\[-.5ex] Michigan State University, East Lansing, Michigan, USA}

\begin{document}

\maketitle


\begin{abstract}
Iterative deblurring, notably the Richardson-Lucy algorithm with and without regularization, is analyzed in the context of nuclear and high-energy physics applications.  In these applications, probability distributions may be discretized into a few bins, measurement statistics can be high, and instrument performance can be well understood.  In such circumstances, it is essential to understand the deblurring first without any explicit noise considerations.  We employ singular value decomposition for the blurring matrix in a low-count pixel system.  A strong blurring may yield a null space for the blurring matrix.  Yet, a nonnegativity constraint for images built into the deblurring may help restore null-space content in a high-contrast image with zero or low intensity for a sufficient number of pixels.  For low-contrast images, control over null-space content can be achieved through regularization.  When regularization is applied, the blurred image is, in practice, restored to one that is still blurred but less than the starting image.
\end{abstract}

\section{Introduction}
Blurring is common in image acquisition, making the images less sharp and clear.   When the blurring process is understood, it may be possible to remove, fully or partially, the blurring effects or to deblur the images.  Here, we will analyze what a standard iterative deblurring methodology can and cannot do, with an eye toward improving the quantitative results of particle yield measurements in nuclear and high-energy physics.  In the latter areas, the deblurring gained attention relatively recently \cite{dagostini_multidimensional_1995, vargas_unfolding_2013, danielewicz_deblurring_2022, nzabahimana_deconvoluting_2023, xu_reconstruction_2024}.  The measurement statistics can be higher than in many optical applications, making non-nominal noise analysis relevant.  At the same time, the number of bins in these areas, equivalent to pixels, can be modest.

 A convenient framework in which the deblurring can be understood is that of the Singular Value Decomposition (SVD) for the blurring matrix, already used in the past in the literature \cite{hansen_deblurring_2006, reichel_introduction_2013, zuninga_singular_2021}.  With some advantages that the Richardson-Lucy~(RL) \cite{richardson_bayesian-based_1972, lucy_iterative_1974} and related iterative methods \cite{khan_iterative_2013, dey_richardsonlucy_2006, danielewicz_deblurring_2022} have over other iterative methods for the restoration of brightness \cite{biemond_iterative_1990, vankawala_survey_2015, landweber_iteration_1951}, we will concentrate here on the RL method.  One novel aspect of our work will examine the contrast's impact on the objects imaged in the deblurring process.  Further, we will actively seek cases where deblurring fails, understand the failures, and predict when they may occur.  Small or vanishing eigenvalues from SVD will play a role there.  In the present work, we will focus on deblurring when no explicit noise is present in the problem, which is already complex.  In future work, we shall consider practicalities around the noise in the problem, for which the present work will provide a base.

This paper is organized as follows.  Section 2 discusses general blurring concepts and formulates a schematic model for explorations.  We decompose the blurring matrix using SVD.  Section 3 describes selected deblurring methods, including the RL and Landweber (LW) algorithms, and illustrates their practical operation.  Section~4 explores the limits of the RL deblurring with and without regularization.  Section~5 discusses how concurrent nuclear physics efforts connect to this paper's explorations.  We conclude in Section 6.

\section{Blurring Matrices and Their Singular Value Decomposition}
The blurring relation can be stated in the form of the  equation \cite{lucy_iterative_1974} 
 \begin{align}\label{main}
 \mathbf{g}=T \mathbf{F} \, .
 \end{align}
 Here, $ \mathbf{g} \in \mathbb{R}^{n}$ is the vector of values for the blurred image, $T\in \mathbb{R}^{n\times n}$ is the blurring matrix or transfer matrix, and $ \mathbf{F} \in \mathbb{R}^{n}$ is the vector for the blur-free image.  Probability conservation for the blurring matrix implies that the sum of all the elements of each column (index $i$) is precisely equal to~1:
 $$
 \sum_{j=0}^{n-1} T_{j i}=1 \, .
 $$

Without practical loss of generality, we shall consider that the indices correspond to locations in one dimension, representing pixels or other distribution bins.  We shall adopt cyclic boundary conditions within that space, where the vector and matrix indices are computed modulo~$n$.  To fix the attention, we shall consider a couple of exemplary blurring (or point-spreading) functions, three- and five-bin, out of which we shall consider the five-bin one most often:\\[-.5ex]

\noindent
Three-bin
 \begin{align}\label{eq:t3}
    T_{ji}=
    \begin{cases}
    \begin{array}{cl}
    0.550 \, , & \text{for } j = i \, , \\
    0.225 \, , & \text{for } j = (i \pm 1) \mod n \, ,\\
    0 \, , & \text{otherwise.}\\
    \end{array}
    \end{cases}
 \end{align}
 
\begin{figure}[!hptb]
    \centering
    \begin{subfigure}[b]{0.46\linewidth}        
        \centering
        \includegraphics[width=\linewidth]{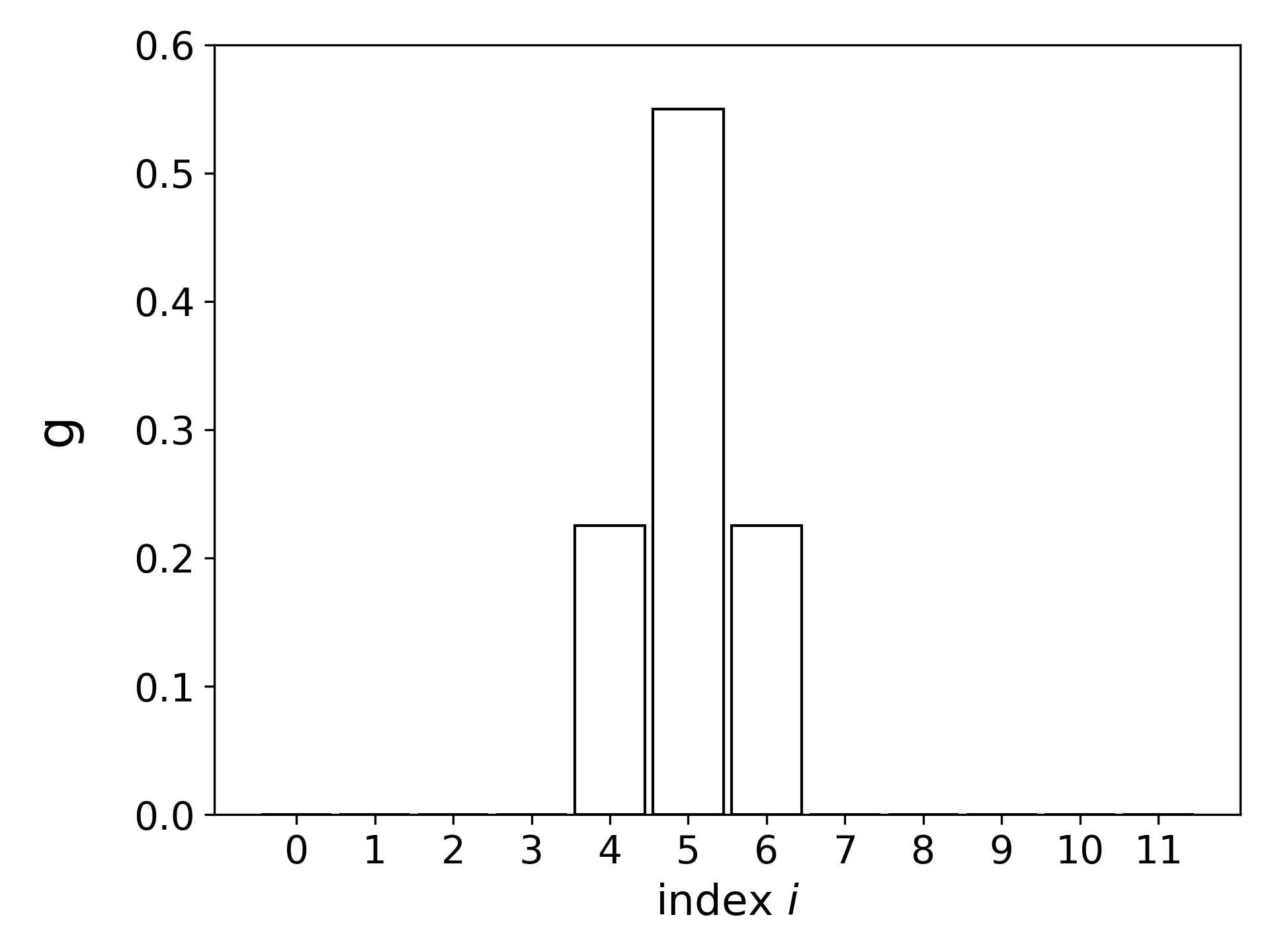}
        \caption*{(a)}
    \end{subfigure}
    \hspace*{1em}
    \begin{subfigure}[b]{0.46\linewidth}        
        \centering
        \includegraphics[width=\linewidth]{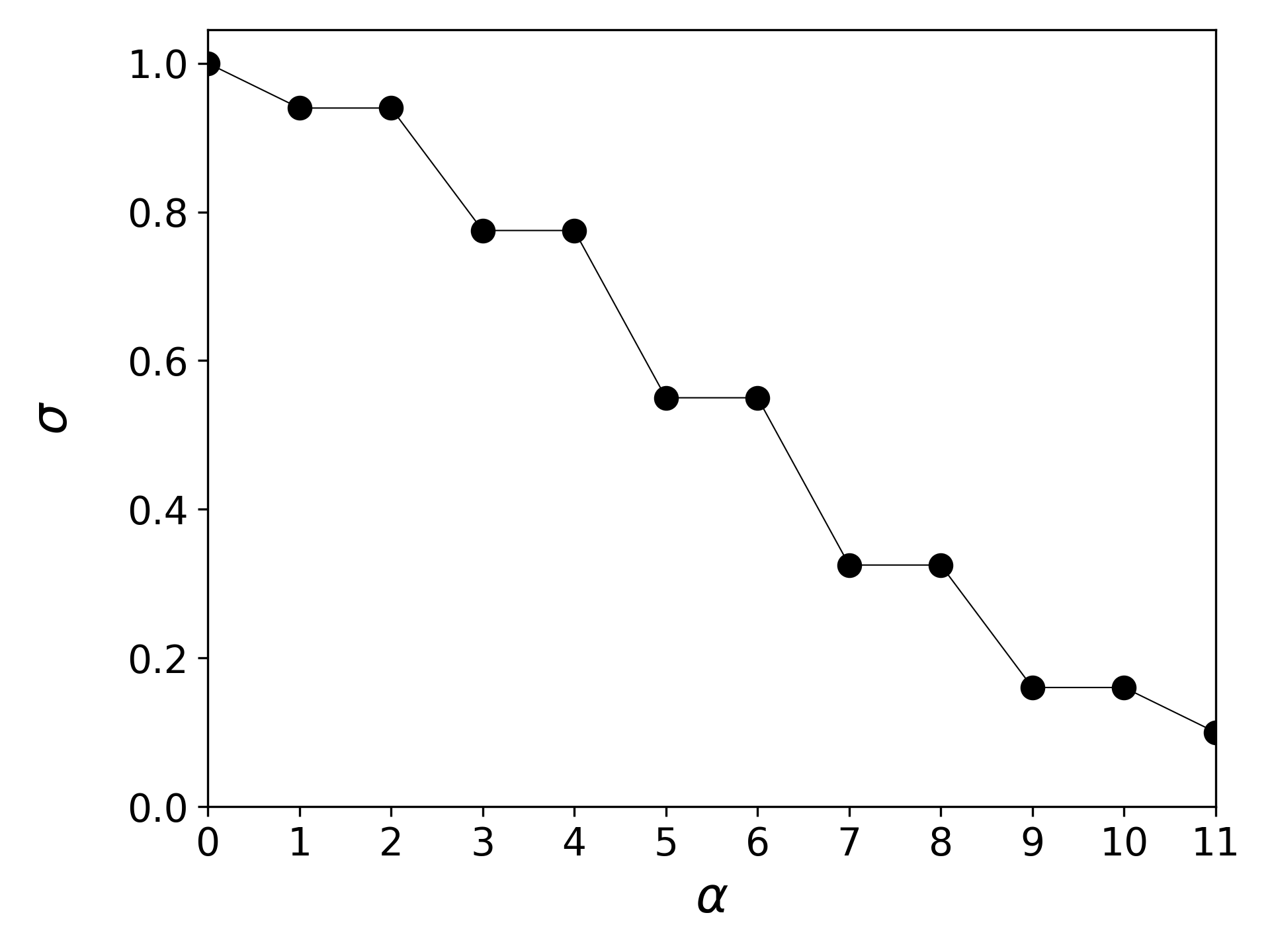}
        \caption*{(b)}
    \end{subfigure}
    \caption{Illustration of the three-bin blurring function.  (a) Image from blurring an original with pixel 5 at 1 and other at 0, for $n=12$.  (b) Singular values of the blurring function, ordered by their magnitude.    Lines join the values to guide the eye. }  
        \label{fig:threebin}
 \end{figure}

\begin{figure}[!hptb]
    \centering
    \begin{subfigure}[b]{0.46\linewidth}        
        \centering
        \includegraphics[width=\linewidth]{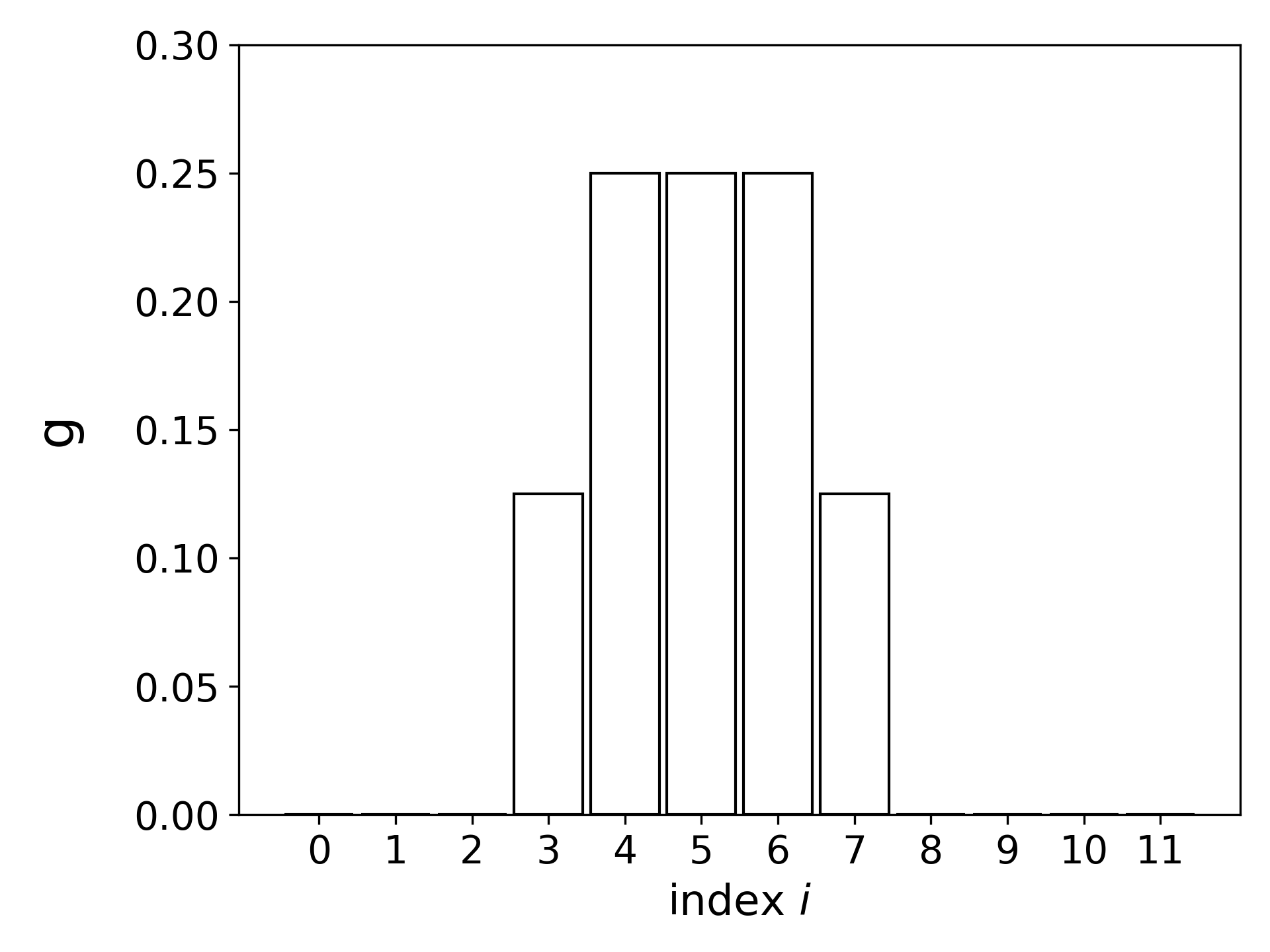}
        \caption*{(a)}
    \end{subfigure}
     \hspace*{1em}
    \begin{subfigure}[b]{0.46\linewidth}        
        \centering
        \includegraphics[width=\linewidth]{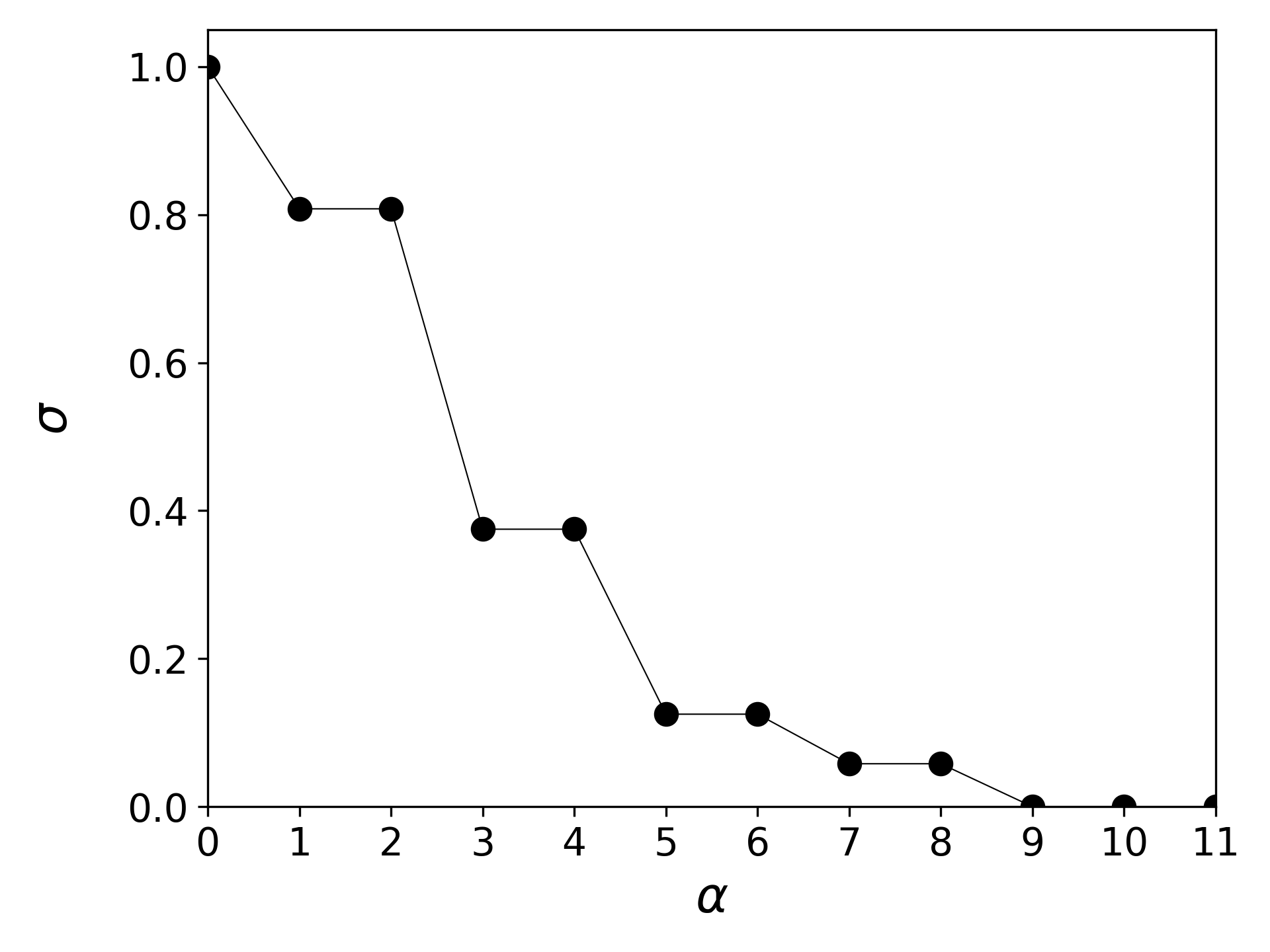}
        \caption*{(b)}
    \end{subfigure}
    \caption{Illustration of the five-bin blurring function.  (a) Image from blurring an original with pixel 5 at 1 and other at 0, for $n=12$.  (b) Singular values of the blurring function, ordered by their magnitude.  Lines join the values to guide the eye.}
    \label{fig:fivebin}
 \end{figure}

\noindent  
Five-bin
 \begin{align}\label{eq:t5}
    T_{ji}=
    \begin{cases}
    \begin{array}{cl}
    0.250 \, , & \text{for } j=i \, , \\
    0.250 \, , & \text{for } j=(i \pm 1) \mod n \, , \\
    0.125 \, , & \text{for } j= (i \pm 2) \mod n \, , \\
    0 \, , & \text{otherwise.}\\   
    \end{array}
    \end{cases}
 \end{align}
In the grayscale convention,  the pixel brightness values span the interval [0,1], with~1 representing white and 0 representing black.  In deblurring, it is important not to end up with values outside of that interval, as these would lack physical sense.  We illustrate the three- and five-bin blurring functions in Fig.~\ref{fig:threebin}(a) and Fig.~\ref{fig:fivebin}(a), respectively, where we show the action of these functions on an original image where only pixel~5 has a value of 1, and the others have 0.  The overall size of the space is $n=12$.

Upon applying the transposed blurring matrix to both sides of Eq.~\eqref{main}, we find the norm equation
\begin{equation}\label{g}
    T^\top  \mathbf{g} = T^\top  T \, \mathbf{F} \equiv M \mathbf{F} \, .
\end{equation}
If we apply SVD to $T$, we get
\begin{equation}
    T= \sum_{\alpha=0}^{n-1} \sigma_\alpha \, \mathbf{u}_\alpha \, \mathbf{v}_\alpha^\top  \, ,
\end{equation}
where $\mathbf{u}$ and $\mathbf{v}$ are, respectively, the separately orthonormal left and right singular vectors, and we adopt the convention of the singular values $\sigma$ to be nonnegative and ordered in nonincreasing order.  Then for the positive definite matrix $M$ we arrive at the following spectral decomposition
\begin{equation}\label{3}
    M = T^\top \, T  =\sum_{\alpha, \beta} \sigma_\alpha \, \sigma_\beta \, \mathbf{v}_{\alpha} \, \mathbf{u}_{\alpha}^\top \, \mathbf{u}_{\beta} \, \mathbf{v}_{\beta}^\top 
 =\sum_\alpha \sigma_\alpha^2 \, \mathbf{v}_{\alpha} \, \mathbf{v}_{\alpha}^\top \, .
\end{equation}

Generally, for original and deblurred images, we shall employ a decomposition in terms of the right singular vectors of the blurring matrix and for the blurred images - in terms of the left:
\begin{align}\label{al}
    \begin{cases}
  \mathbf{ F}=\sum_{\alpha}^{} F_{\alpha} \, \mathbf{v}_{\alpha} \, , \\ 
    \mathbf{g}=\sum_{\alpha}^{} g_{\alpha} \, \mathbf{u}_{\alpha} \, .
    \end{cases}
\end{align}
With \eqref{main}, we then find
\begin{equation}
    \mathbf{g} = T \mathbf{F} = \sum_{\alpha} \sigma_\alpha \, F_\alpha \, \mathbf{u}_\alpha \, ,
\end{equation}
or $g_\alpha = \sigma_\alpha \, F_\alpha$.

The three-  and five-bin blurring matrices, Eqs.~\eqref{eq:t3} and \eqref{eq:t5}, which we use here as examples, have some peculiar features.  These features are often encountered in practice, possibly as approximate rather than exact features.  Thus, these blurring matrices are symmetric, meaning their singular vectors are matrix eigenvectors.  Further, with the matrix elements of~$T$ depending only on the absolute difference of indices, that dependence extends to the corresponding~$M$ matrices.  This implies two types of invariance for the matrices~$ T$ and~$ M$, which affect the singular vectors.  One is the invariance under a shift of matrix indices by some integer.  With this, a singular vector, corresponding to some singular value, remains a singular vector, corresponding to the same singular value, after a shift in its indices by an integer.  The second invariance is under a reflection of the indices around any chosen index.  With this, a singular vector, corresponding to some singular value, remains a singular vector, corresponding to the same singular value, after a reflection of its indices around any chosen index.  With these two symmetries, any singular value is at least doubly degenerate unless the corresponding singular vector transforms onto itself, up to a factor of -1, under these transformations.
 
\begin{figure}[!hptb]
    \centering
    
    \begin{subfigure}{0.46\textwidth}
        \includegraphics[width=\linewidth]{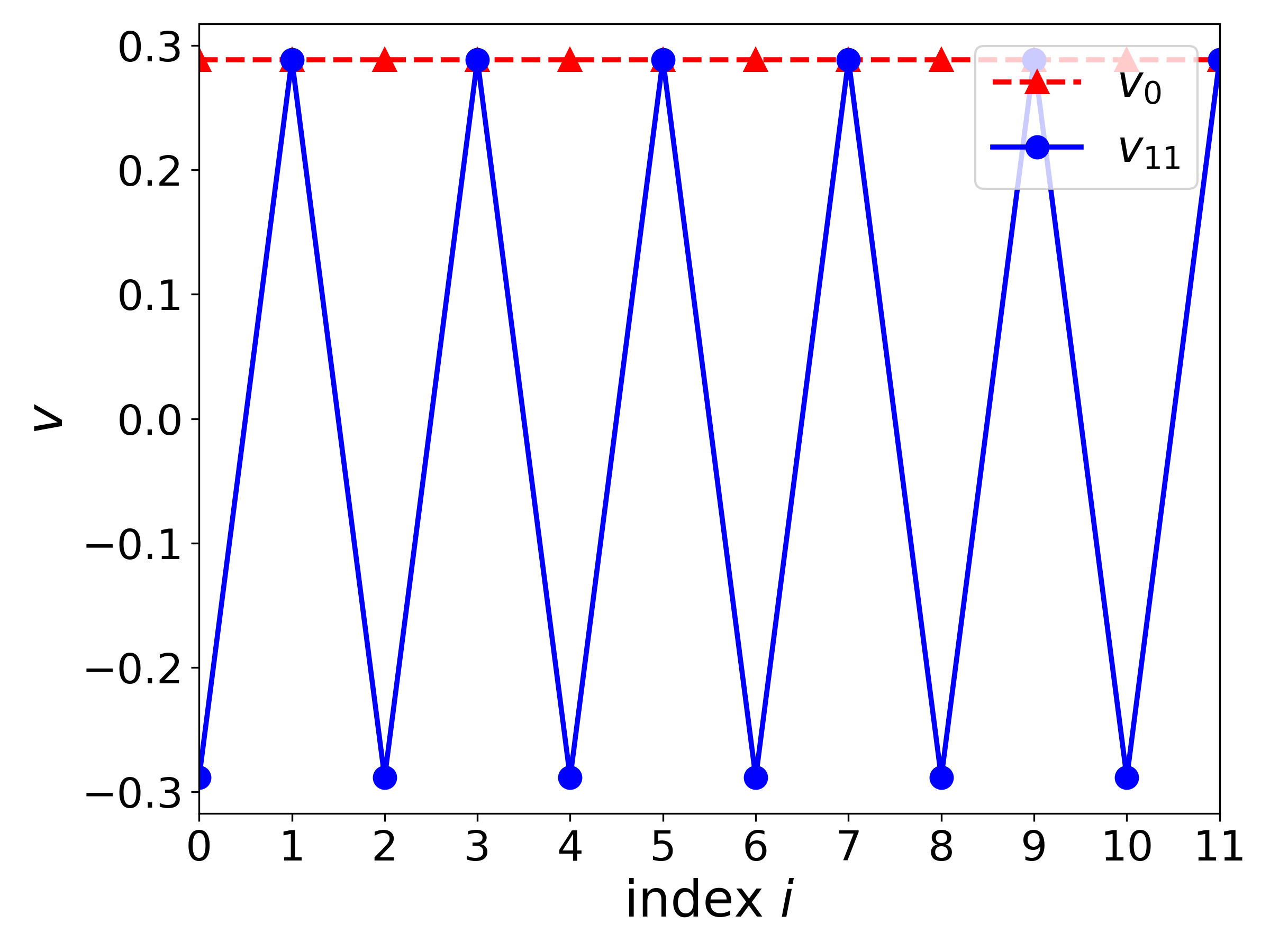}
        \caption*{(a)}
    \end{subfigure}
     \hspace*{1em}
    \begin{subfigure}{0.46\textwidth}
        \includegraphics[width=\linewidth]{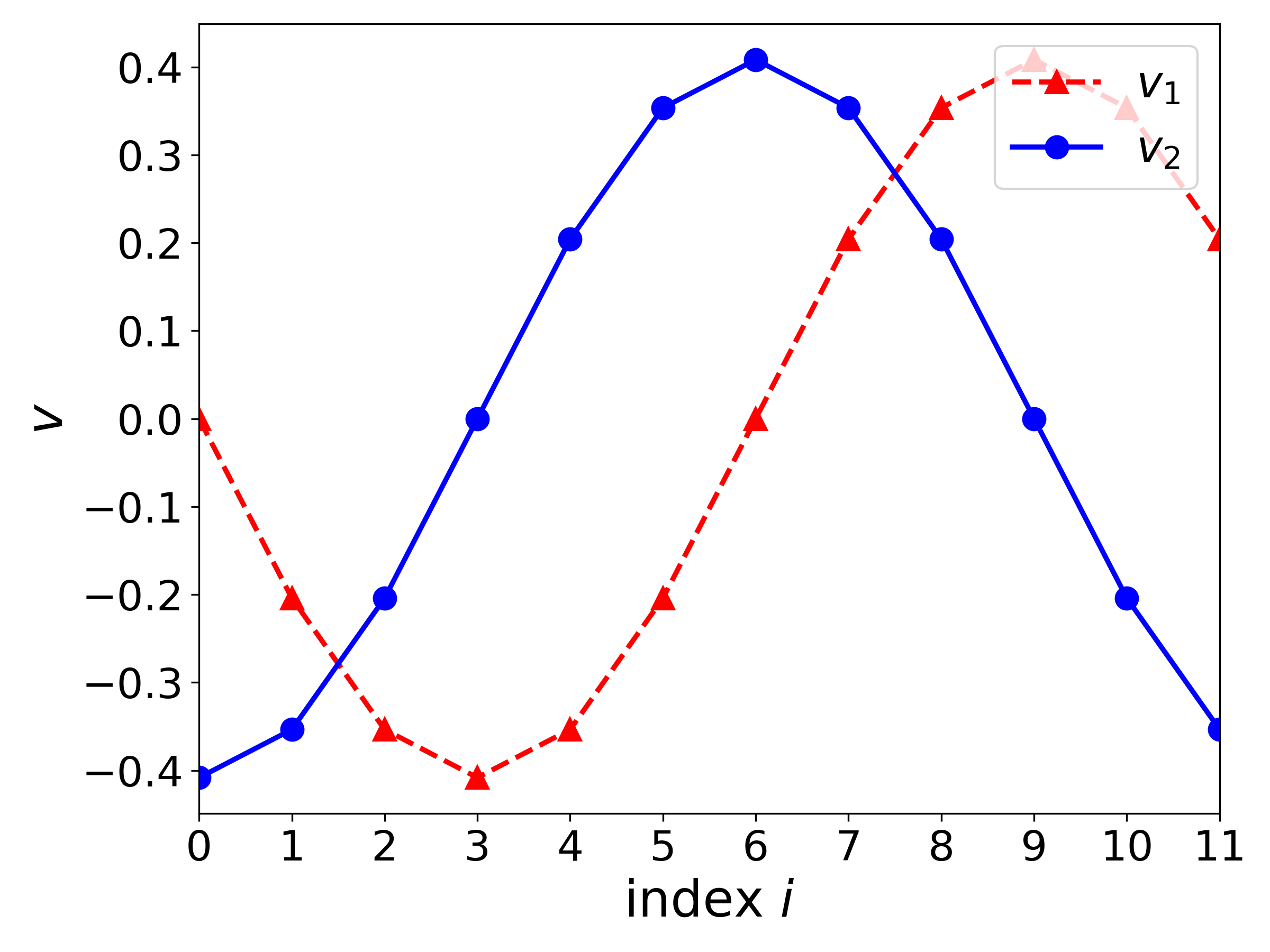}
         \caption*{(b)}
    \end{subfigure}
    
    \vspace{0.6em} 
    
    \begin{subfigure}{0.46\textwidth}
        \includegraphics[width=\linewidth]{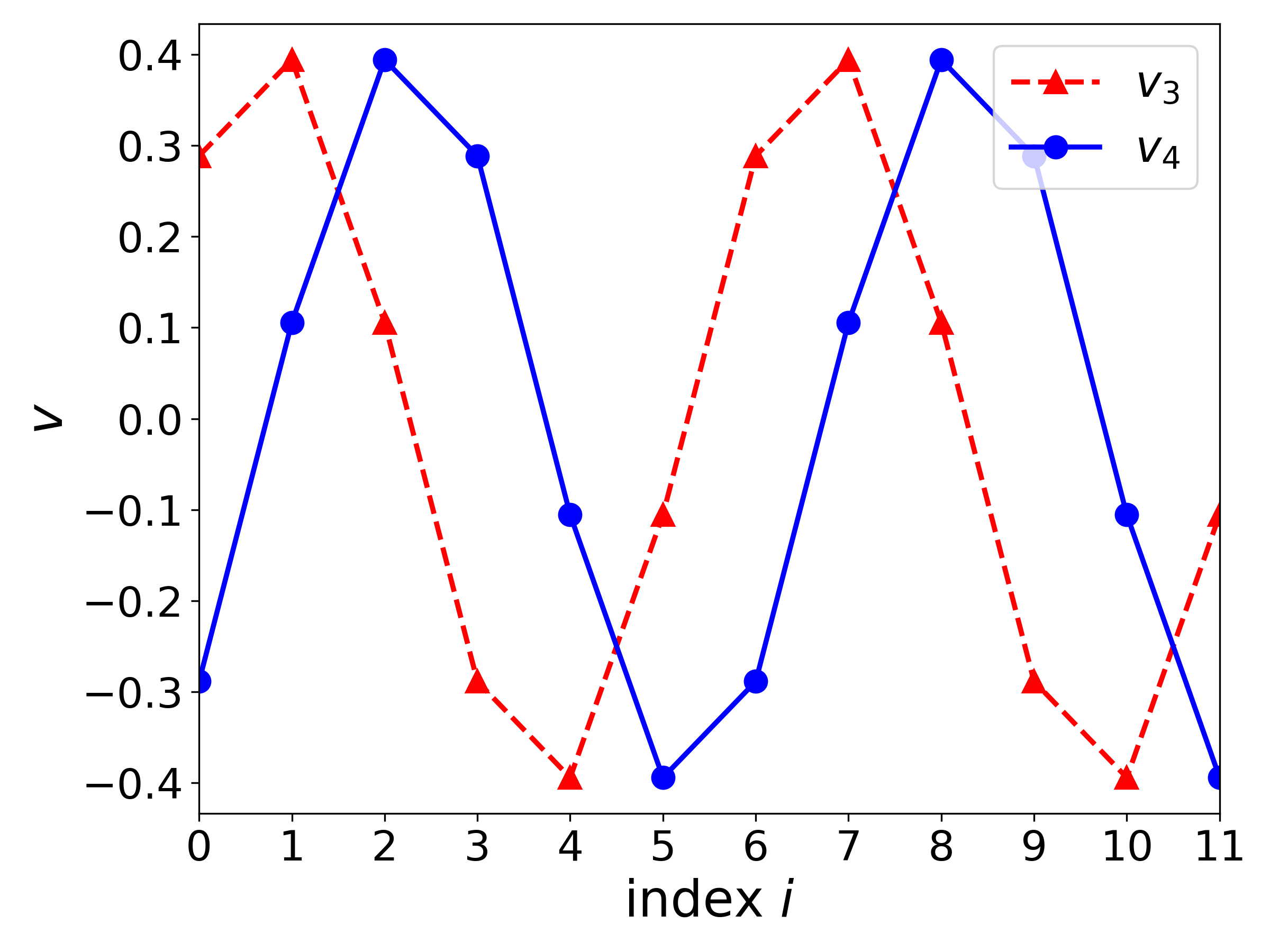}
         \caption*{(c)}
    \end{subfigure}
    \hspace*{1em}
    \begin{subfigure}{0.46\textwidth}
        \includegraphics[width=\linewidth]{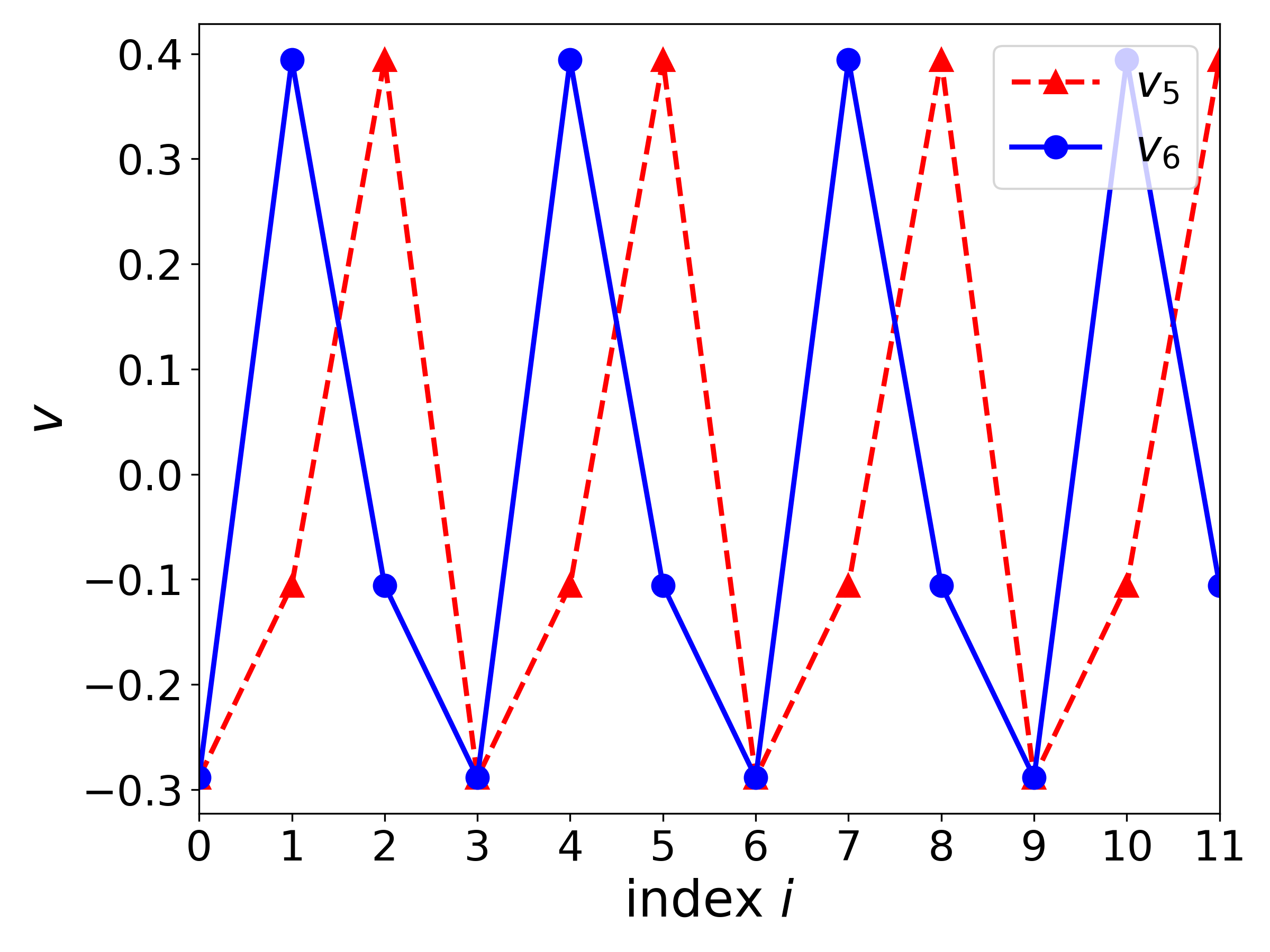}
         \caption*{(d)}
    \end{subfigure}
    
    \vspace{0.4em} 
    
    \begin{subfigure}{0.46\textwidth}
        \includegraphics[width=\linewidth]{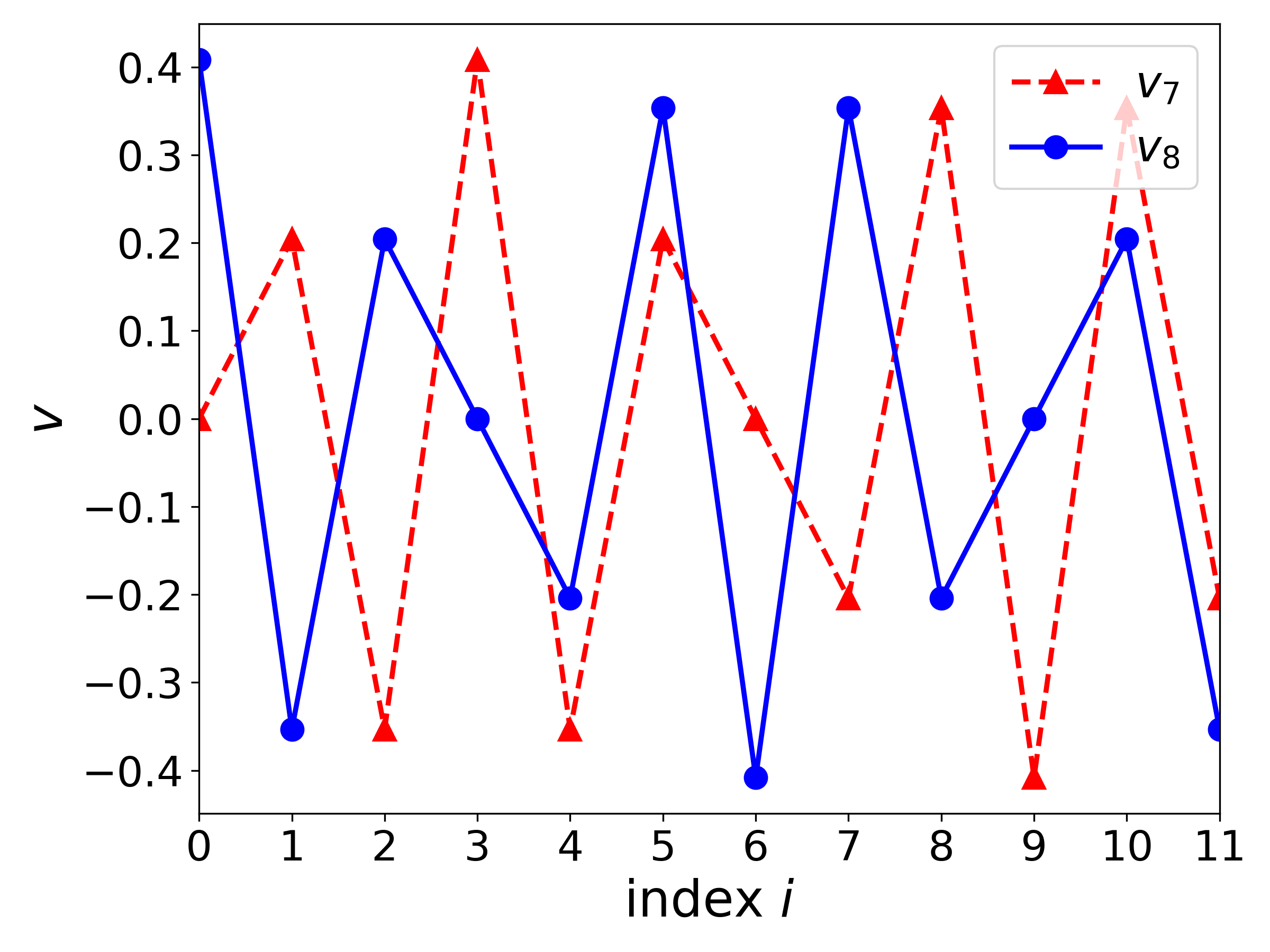}
         \caption*{(e)}
    \end{subfigure}
     \hspace*{1em}
    \begin{subfigure}{0.46\textwidth}
        \includegraphics[width=\linewidth]{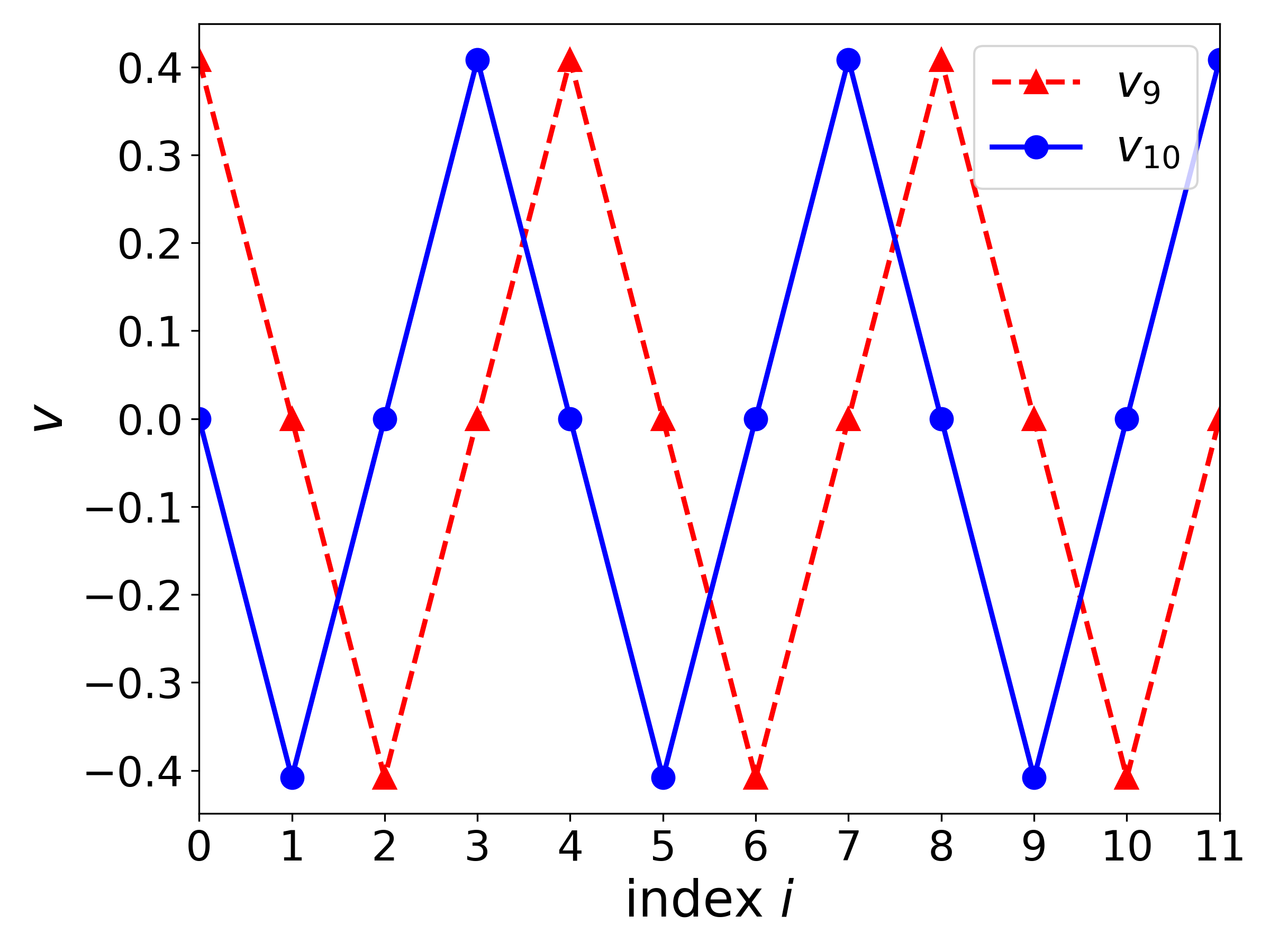}
        \caption*{(f)}
    \end{subfigure}
    \caption{Singular vectors for the five-bin blurring function.  Within degenerate spaces, pairs of vectors are chosen to transform onto the other or itself under reflection and/or translation.  Panel (a) shows the two unique vectors that transform onto themselves under the transformations.}    
       \label{fig:singulars}
\end{figure}

The spectra of singular values for the two exemplary blurring matrices are shown in Figs.~\ref{fig:threebin}(b) and \ref{fig:fivebin}(b), and the singular vectors for the five-bin function are shown in Fig.~\ref{fig:singulars}.  Double degeneracy is seen for most of the singular values.  For Fig.~\ref{fig:singulars}, we choose the corresponding pairs of singular vectors in the degenerate subspaces so that some symmetry behind the degeneracies is evident, either a reflection or translational one, or both.  In each of the spectra, in Figs.~\ref{fig:threebin}(b) and \ref{fig:fivebin}(b), a nondegenerate singular value of 1 is seen.  This value corresponds to the vector that is constant across the pixels.  Under blurring, this vector transforms into itself.  Also, this vector is invariant under the symmetries.  The spectrum of singular values falls off faster for the five- than for the three-bin function, demonstrating a more significant deterioration of the information under the blurring with a broader function.  In Fig.~\ref{fig:singulars}, it is evident that a faster variation of the vectors over pixels is generally associated with lower singular values.  The singular vector with the most rapid variation, alternating between the same positive and negative values from one pixel to another, is invariant under the symmetries.  It is the other one that may be nondegenerate.  Also, under blurring, it transforms into itself, up to a factor, because it is uniquely alternating and reflection-symmetric about every pixel, which will be preserved by any blurring that depends only on the difference of indices.  With a broadening of the blurring function, this vector will likely be the first to be nullified.  In Fig.~\ref{fig:threebin}(b), the singular value for that vector, with index 11, approaches zero, and the vector is nullified when the blurring matrix evolves into $T_{jj} = 0.5$ and $T_{ji}=0.25$, for $j = (i\pm 1) \mod 12$.  In Fig.~\ref{fig:fivebin}(b), the maximally alternating vector is accompanied in the null space by two more vectors, cf., panels (a) and (f) in Fig.~\ref{fig:singulars}.


\section{Deblurring Methods}
\begin{figure}[!htb]
    \centering 
    \begin{subfigure}[b]{0.45\linewidth}        
        \centering
        \includegraphics[width=\linewidth]{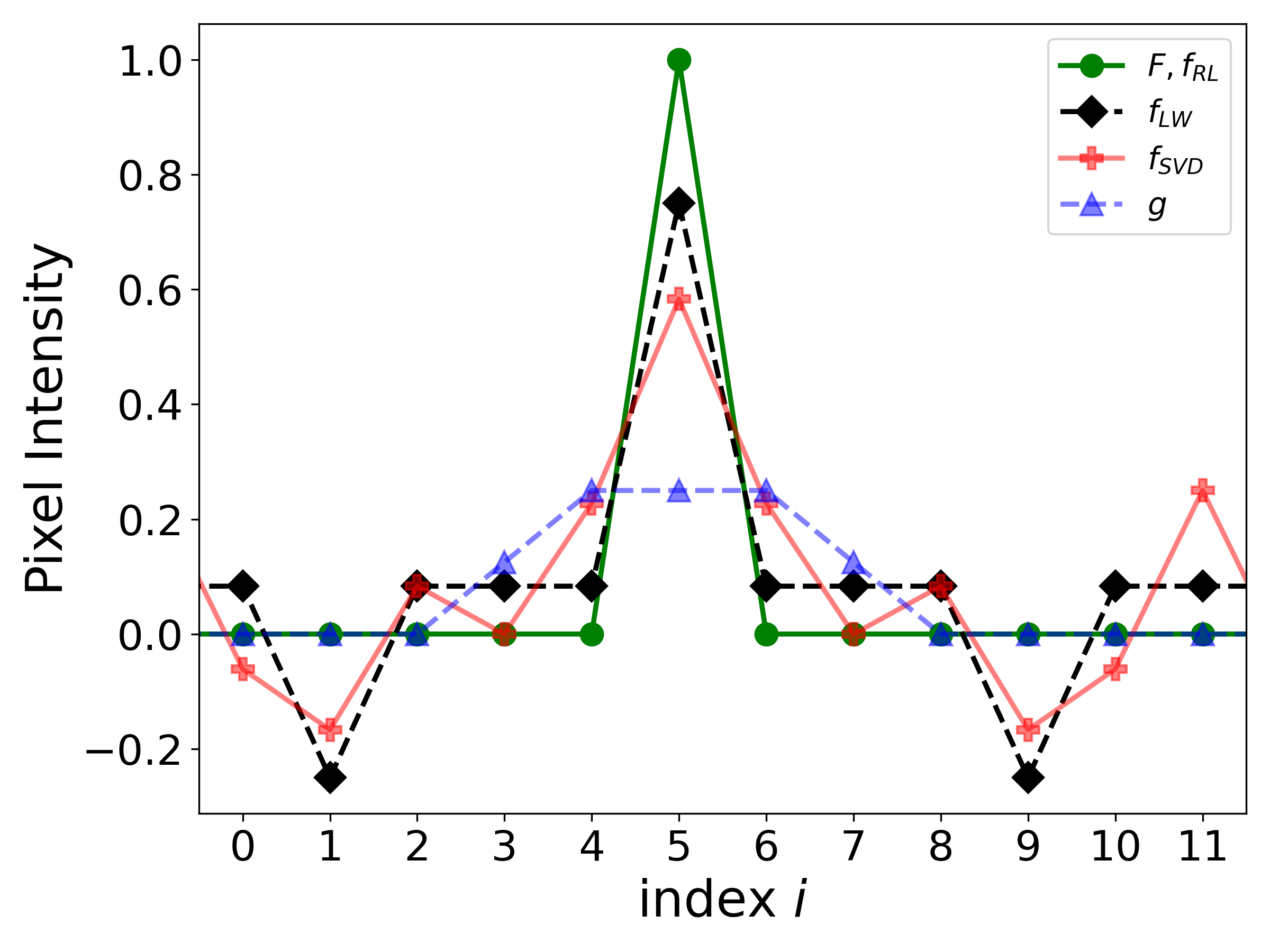}\\
        \caption*{(a)}
    \end{subfigure}
     \hspace*{1em}
    \begin{subfigure}[b]{0.47\linewidth}        
        \centering
    \includegraphics[width=\linewidth]{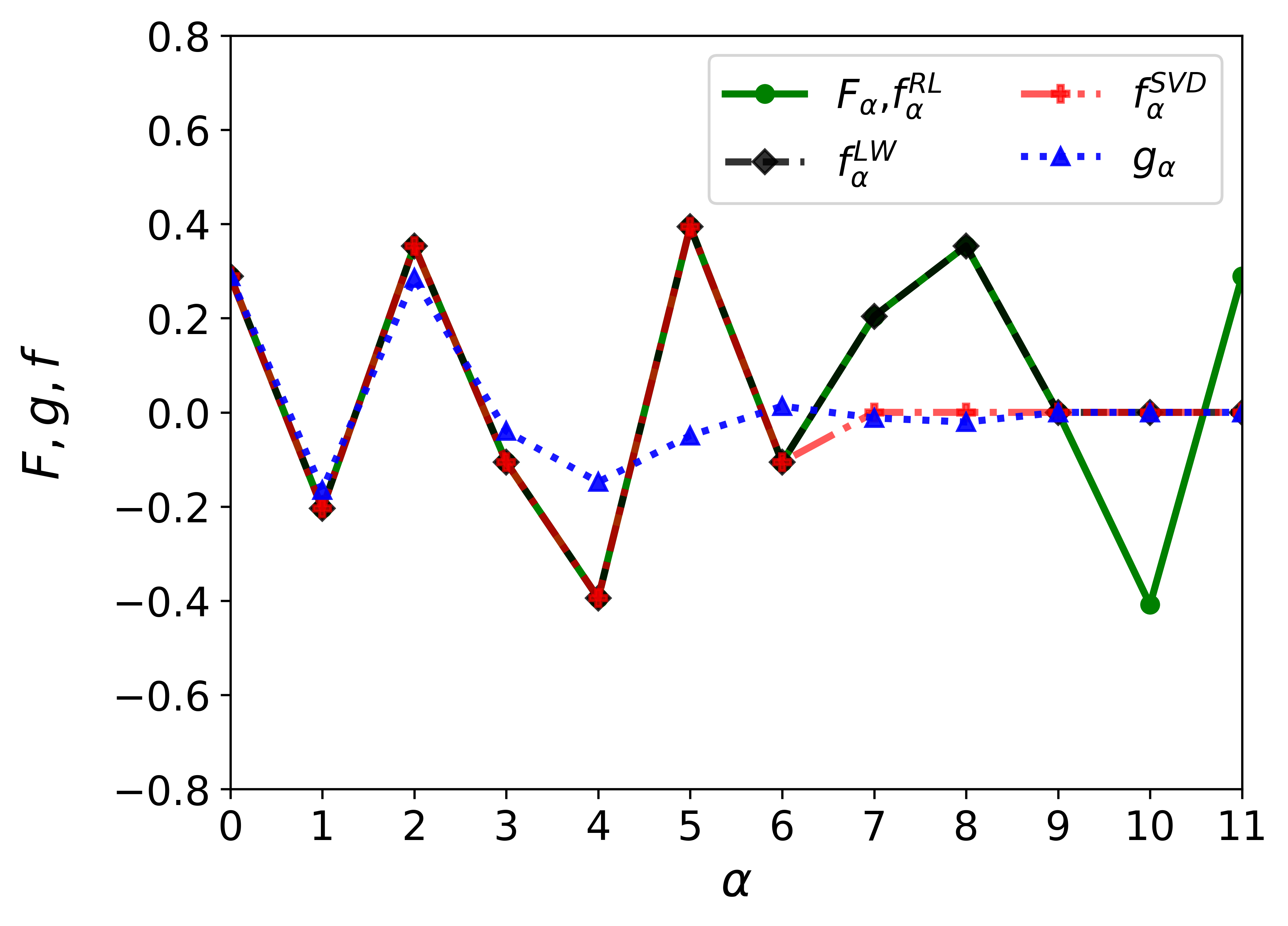}
        \caption*{(b)}
    \end{subfigure}
    \caption{Deblurring test with a single pixel with index 5 at the intensity 1 and all other pixels at the intensity of 0.   (a) The original image $F$ (circles), image $g$ (triangles) from blurring $F$ with the 5-pixel blurring function, and the restored images $f$ from different deblurring methods are displayed.  As the image restored by the RL method is indistinguishable from the original, the circles represent that image, too.  The image restored using the LW method is shown as diamonds.  Finally, the image restored with SVD, where the operational space was reduced ($m=6$ in Eqs.~\eqref{eq:SVDinverse} and \eqref{eq:SVD0}) compared to LW, is shown as circles.  The symbols corresponding to each case are connected with lines. (b) Coefficients of decomposition in the basis of singular vectors vs singular value index for the original and blurred images and images restored with the RL and LW methods.}
    \label{fig:delta} 
 \end{figure}

We next turn to deblurring.  We illustrate the deblurring in Fig.~\ref{fig:delta}(a).  We take a $\delta$ distribution where only one bin gets intensity 1, and the others get 0.  We blur the original image with the 5-bin function \eqref{eq:t5} and deblur the result using the RL method \cite{richardson_bayesian-based_1972, lucy_iterative_1974}, details of which will be discussed.  The blurred and restored images are shown in Fig.~\ref{fig:delta}(a) together with the original.  The eye cannot distinguish between the RL-restored and original images, so they are represented with a single set of graphics.  In addition, the intensity values for the images are shown in Table~\ref{table:delta}, and within the precision of the machine we employ, they cannot be distinguished.  This may be perceived as astounding, as the singular vectors with the fastest bin-to-bin variation span the null space for the blurring matrix, and they most certainly contribute to the SVD of the $\delta$ distribution.

\begin{table}[!htb]
    \centering 
  \begin{tabular}{|c|c|c|c|c|c|}
    \hline
    {$i$} &  $F_{i}$ &      $g_{i}$ &           $f_{i}^{RL}$ & $f_{i}^{LW}$ &     $f_{i}^{SVD}$       \\
    \hline
0  &  0 &  0.000 &   0.000 &  0.083 & -0.061 \\
1  &  0 &  0.000 &   0.000 & -0.250 & -0.167 \\
2  &  0 &  0.000 &   0.000 &  0.083 &  0.083 \\
3  &  0 &  0.125 &   0.000 &  0.083 &  0.000 \\
4  &  0 &  0.250 &   0.000 &  0.083 &  0.228 \\
5  &  1 &  0.250 &   1.000 &  0.750 &  0.583 \\
6  &  0 &  0.250 &   0.000 &  0.083 &  0.228 \\
7  &  0 &  0.125 &   0.000 &  0.083 &  0.000 \\
8  &  0 &  0.000 &   0.000 &  0.083 &  0.083 \\
9  &  0 &  0.000 &   0.000 & -0.250 & -0.167 \\
10 &  0 &  0.000 &   0.000 &  0.083 & -0.061 \\
11 &  0 &  0.000 &   0.000 &  0.083 &  0.250 \\
   \hline
    \end{tabular}
    \caption{Bin intensities in deblurring the $\delta$-distribution centered at index 5, illustrated in Fig.~\ref{fig:delta}(a).  From left to right, the columns represent the bin index~($i$), and then the intensity for the original image ($F_i$), for the 5-bin blurred image ($g_i$), and the image restored with three different methods: RL ($f_i^{RL}$),  LW ($f_i^{LW}$)  and  SVD ($f_i^{SVD}$).}
    \label{table:delta}
\end{table}

In the following, we review exemplary deblurring methods and their potential utility for quantitative research analyses.  After this section, we shall examine under what circumstances the deblurring can yield faithful results, such as in the example above, and we shall seek conditions and a manner under which the deblurring can fail.  Iterative methods aim to restore a blurred image to its original form by iteratively improving estimates of the blur-free image.  These methods take as input an initial guess of the blur-free image, the blurring matrix, and a criterion for the number of iterations.  Typically, the initial guesses are the uniform or blurred image.  Since the blur is to be suppressed at each iteration, naively, a more accurate output is expected with more iterations—in practice, this may not be the case. 

While discussing the iterative methods below, we use $k$ as an iteration index.


\subsection{Richardson-Lucy method, classical and regularized}
The classical RL algorithm iteratively upgrades the restored image $f$, given the blurred image $g$, according to \cite{richardson_bayesian-based_1972, lucy_iterative_1974, danielewicz_deblurring_2022}
\begin{equation}\label{eq:RLiteration}
f_i^{(k+1)}=A_i^{(k)} \, f_i^{(k)} \,.
\end{equation}
Here, the amplification factor $A_i$ for the $i$'th pixel intensity is
\begin{equation}\label{eq:RLamplification}
    A_i^{(k)} = \frac{1}{q_i} \sum_{j} T_{ji} \frac{g_j}{g_j^{(k)}} \, ,
\end{equation}
$i,j, m =0,1,2, \ldots, n-1$, $g_j^{(k)} = \sum_{m} T_{jm} f_m^{(k)}$ is the blurred-image prediction for the restored image $f^{(k)}$ at the $k$'th step of restoration, $T$ is the before introduced point spread function or blurring matrix, and $q_i=\sum_{j} T_{j i}$ serves the normalization.  The normalization factors for the introduced blurring-function examples are unity, $q_i \equiv 1$.

The regularized RL algorithm modifies the iteration equation to
\begin{align}\label{tt}
   f_i^{(k+1)}=A_i^{(k)} \, I^{k}_{i} \, f_i^{(k)} \quad i,j, m =0,1,2, \ldots, n-1 \, ,
\end{align}
where the regularization factor is \cite{dey_richardsonlucy_2006, danielewicz_deblurring_2022}
\begin{equation}
\begin{split}
    I^{(k)}_{i} & = \frac{1}{1-\Lambda  \div \left(\frac{\nabla f_i^{(k)}}{|\nabla f_i^{(k)}|}\right)}\\ & =
    \begin{cases}
     \frac{1}{1-\Lambda} \, ,& \text{ if }  f^{(k)}_{i} < f^{(k)}_{[(i-1) \, \mathrm{mod} \, n] , \, [(i+1) \, \mathrm{mod} \, n] } \, , \\
     \frac{1}{1+\Lambda} \, ,& \text{ if }  f^{(k)}_{i} > f^{(k)}_{[(i-1) \, \mathrm{mod} \, n] , \, [(i+1)  \mathrm{mod} \, n] } \, ,  \\
     1 \, , &\text{  }  \text{ otherwise } \, .
    \end{cases}   
\end{split}  
    \label{eq:I=}
    \end{equation}
for $i=0,1,2, \ldots, n-1$.  Here, the r.h.s.\ is the one-dimensional version \cite{danielewicz_deblurring_2022} of the center expression \cite{dey_richardsonlucy_2006} in which $\nabla$ is the gradient and $|\nabla f_i^{(k)}|$ is $\ell_1$-norm; $\Lambda$~is a small positive regularization coefficient.  The role of the regularization factor is to suppress an uncontrolled rise, during iterations, in the contributions to the restored image from the fastest-varying singular vectors corresponding to vanishing or minimal singular values.  The typical signature of such an uncontrolled rise is developing a bin-to-bin seesaw pattern in the restored image.  As SVD might not be applied in parallel to the RL-restoration and the direct connection of the regularization factor $I$ to SVD is not obvious, it is common to make $\Lambda$ small enough that no pathological signs in the restored image appear.

\begin{figure}[!htb]
   \centering 
   \includegraphics[width=.46\textwidth]{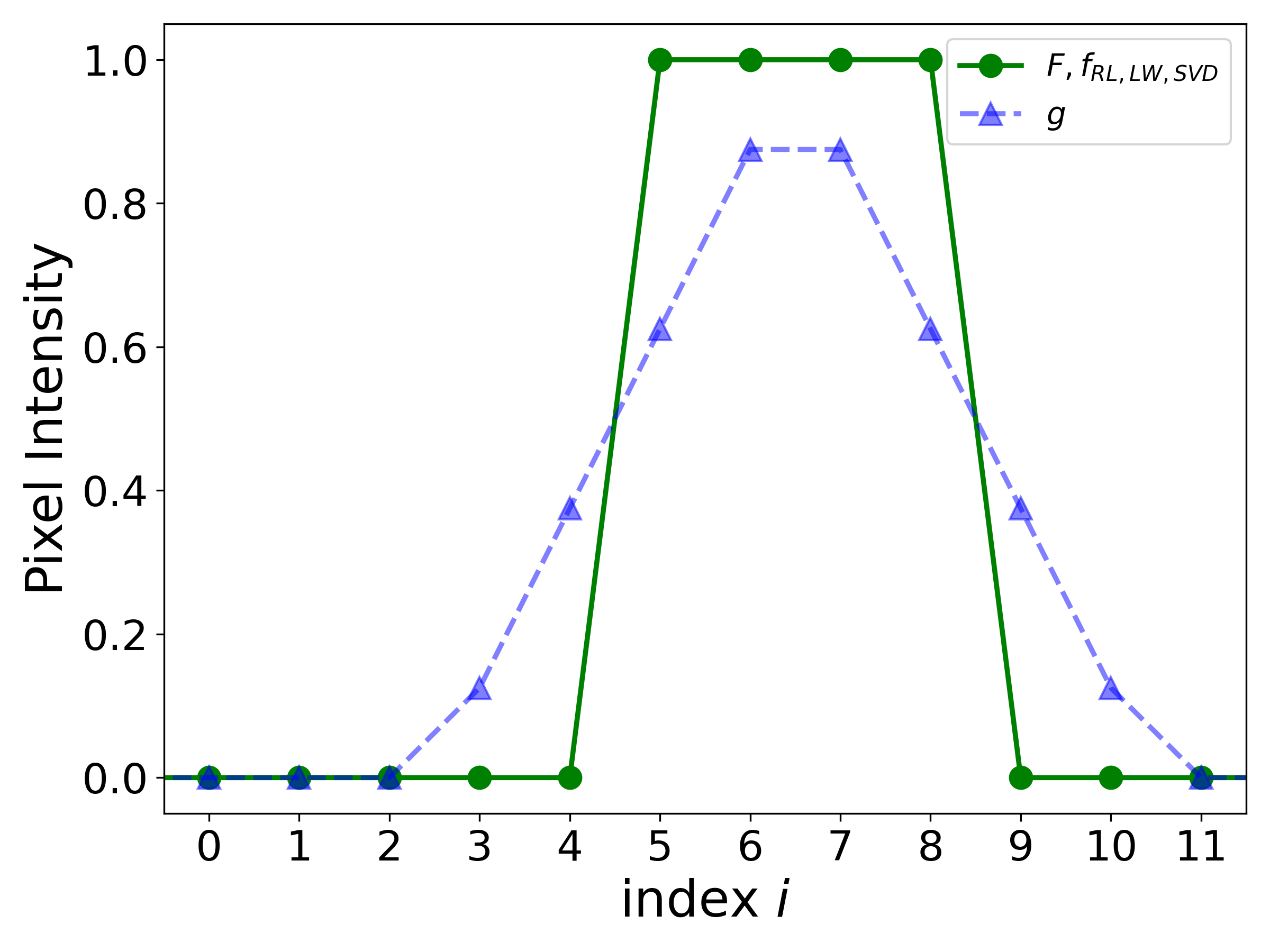}
   \caption{Deblurring test with a 4-pixel ramp.  In the original image, the intensity is 1 for the pixels 5-8 and 0 for the others.  The operational space for SVD is here the full row space of $T$.}
   \label{fig:4PixelRamp} 
\end{figure}

The restoration of the $\delta$-distribution of Fig.~\ref{fig:delta}(a) with the RL method, from the blurred image there, with or without regularization, yields a result that is not distinguishable by eye from the original, so they share the graphical representation in the figure.  The numerical values for the restoration are provided in Table \ref{table:delta}.  Another test to which we subject the deblurring methods is presented in Fig.~\ref{fig:4PixelRamp}.  The original image $F$ is a ramp with pixels 5-8 at intensity 1 and other pixels at intensity 0.  The blurring, yielding $g$, is performed with the 5-bin blurring function.  The results of RL deblurring, without or with regularization, are not distinguishable by eye from the original, but so are the results of deblurring with other methods that we will discuss next.  The numerical values from the ramp tests are shown in Table \ref{table:4PixelRamp}.


\begin{table}[!htb]
    \centering 
  \begin{tabular}{|c|c|c|c|c|c|}
    \hline
    {$i$} &  $F_{i}$ &      $g_{i}$ &           $f_{i}^{RL}$ &  $f_{i}^{LW}$  &    $f_{i}^{SVD}$     \\
    \hline
0  &  0 &  0.000 &   0.000 & 0.000 & 0.000 \\
1  &  0 &  0.000 &   0.000 & 0.000 &0.000 \\
2  &  0 &  0.000 &    0.000 &0.000 & 0.000 \\
3  &  0 &  0.125 &   0.000 &  0.000&  0.000 \\
4  &  0 &  0.375 &    0.000 & 0.000 & 0.000 \\
5  &  1 &  0.625 &   1.000 & 1.000 & 1.000 \\
6  &  1 &  0.875 &   1.000 & 1.000 & 1.000 \\
7  &  1 &  0.875 &   1.000 & 1.000 &  1.000 \\
8  &  1 &  0.625 &   1.000 & 1.000 &  1.000 \\
9  &  0 &  0.375 &  0.000 & 0.000 & 0.000 \\
10 &  0 &  0.125 & 0.000 & 0.000 &  0.000 \\
11 &  0 &  0.000 &  0.000 & 0.000&  0.000 \\
   \hline
    \end{tabular}
    \caption{The numerical values from the 4-pixel ramp deblurring test, with the RL, SVD, and LW methods, illustrated in Fig.~\ref{fig:4PixelRamp}. }
    \label{table:4PixelRamp}
\end{table}

\subsection{Landweber Method}
The Landweber (LW) method \cite{landweber_iteration_1951} can be viewed as the steepest-descent method for minimizing the square deviation between the brightness distribution for the generated and blurred images,
\begin{equation}\label{func:min}
\left\|\mathbf{g}-T \mathbf{f}\right\|^{2}_{2} \, .
\end{equation}
Iterative adjustments in $\mathbf{f}$ take the form
\begin{equation} \label{eq:LW_k1}
    \mathbf{f}^{(k+1)}= \mathbf{f}^{(k)} + \beta T^T (\mathbf{g} - T \mathbf{f}^{(k)}) \, ,
\end{equation}
where $0 < \beta < 2/\|T^{T}T\|$ for convergence. With the $T$ multiplication in \eqref{eq:LW_k1}, it is apparent that the iterations do not change the null-space content in $f$.  The~null-space content may be suppressed completely in $f$ by replacing the minimized quantity~\eqref{func:min} with
\begin{equation}\label{func:min_mod}
\left\|\mathbf{g}-T \mathbf{f}\right\|^{2}_{2} + \epsilon \left\|\mathbf{f} \right\|^{2}_{2} \, ,
\end{equation}
where $0 < \epsilon \ll (\sigma_\alpha)_{min}^2$ and $(\sigma_\alpha)_{min}$ is the smallest singular value of which the contribution to $f$ should not be disturbed.  The iterations are then modified to
\begin{equation}
    \mathbf{f}^{(k+1)}= (1 - \beta \epsilon) \mathbf{f}^{(k)} + \beta T^T (\mathbf{g} - T \mathbf{f}^{(k)}) \, .
\end{equation}

The results from the deblurring of the blurred $\delta$-distribution with the LW method are provided in Fig.~\ref{fig:delta} and Table \ref{table:delta}.  The corresponding results from processing the 4-pixel ramp are provided in Fig.~\ref{fig:4PixelRamp} and Table \ref{table:4PixelRamp}.  While the LW method struggles to restore the $\delta$ distribution, yielding, in particular, negative intensity values, it faithfully restores the 4-pixel ramp.  This may be surprising, and we will return to this issue later in this Section.




\subsection{Deblurring Using Singular Value Decomposition}
The deblurring can be further carried out by employing SVD explicitly. In general, the SVD deblurring~\cite{hansen_deblurring_2006, zuninga_singular_2021, reichel_introduction_2013} seeks the best approximation to the blurred image in a reduced subspace of singular vectors, which in particular should exclude the null space for~$T$, minimizing the square deviation
\begin{equation}\label{func:minm}
\left\|\mathbf{g}-T_m \mathbf{f}\right\|^{2}_{2} \, ,
\end{equation}
where
\begin{equation}
T_m = \sum_{\alpha=0}^{m-1} \sigma_\alpha \, \mathbf{u}_\alpha \, \mathbf{v}_\alpha^\top  \, , 
\end{equation}
and $m \le n$.  The minimization of \eqref{func:minm} is equivalent to the minimization of the square deviation within the reduced subspace only
\begin{equation}\label{funcm:minmm}
\left\|\mathbf{g}_m-T_m \mathbf{f}\right\|^{2}_{2} \, ,
\end{equation}
where
\begin{equation}
\mathbf{g}_m = \sum_{\alpha=0}^{m-1} g_\alpha \, \mathbf{u}_\alpha \, , 
\end{equation}
with the minimum reached for
\begin{equation}
f_\alpha = \frac{g_\alpha}{\sigma_\alpha} \, , \hspace*{2em}\text{for} \hspace*{0.5em} \alpha < m \,.
\label{eq:SVDinverse}
\end{equation}
Outside of the considered subspace, the expansion coefficients for $\mathbf{f}$ are set to zero,
\begin{equation}
f_\alpha = 0   \, , \hspace*{2em}\text{for} \hspace*{0.5em} m \le \alpha < n \,.
\label{eq:SVD0}
\end{equation}

In the tests of this Section, the SVD method yields the same results as LW when its operational space is maximal, i.e., identical to the row space of the blurring matrix~$T$.  Regarding the 4-pixel ramp test with results in Fig.~\ref{fig:4PixelRamp} and Table~\ref{table:4PixelRamp}, we keep the operational space for SVD maximal.  However, in the $\delta$-function test shown in Fig.~\ref{fig:delta} and Table~\ref{table:delta}, we reduce the operational space for SVD to $ m=6$. See Fig.~\ref{fig:singulars} for the corresponding singular vectors.  With the reduction, the restoration's quality significantly deteriorates.  The location and strength of negative intensities in the restored image are similar for $m=8$ (identical to LW within precision) and $m=6$, cf.~Eqs.~\eqref{eq:SVDinverse} and \eqref{eq:SVD0}.

\subsection{Deblurring using Deep Learning}
The above traditional deblurring methods rely on predefined blurring, which can be a limitation when the nature of blurring for the images is not fully understood. Deep Learning approaches can leverage large datasets and Deep Neural Networks (DNNs) to infer the mapping between blurred and original images~\cite{trippe2023let}. Due to its data-driven nature, training the DNNs requires high-quality datasets containing pairs of blurred and blur-free images, often synthesized. The idea is to train DNNs to recover the original image $ \mathbf{F}$ from blurred inputs $ \mathbf{g}$ without explicitly modeling or knowing the blurring function $T$ in Eq.~\eqref{main}.


\subsection{Methods' Assessment}

We call the original images, such as in Figs.~\ref{fig:delta}(a) and~\ref{fig:4PixelRamp}, where a substantial fraction of the pixels has an intensity zero or low compared to the maximal, as having high contrast.  When such images are moderately blurred, the deblurring methods with the nonnegativity constraint for intensity built into the deblurring process, such as RL, generally perform far better than methods that lack the constraint, such as LW or SVD, see the case of the $\delta$-distribution of Fig.~\ref{fig:delta}(a).  In this case, the two latter methods yield negative intensity values in the restored images.  Amazingly, while the $\delta$-distribution contains contributions from the $v_{10}$ and $v_{11}$ vectors in the null space, the nonnegativity constraint allows the RL method to restore their contributions faithfully despite their singular values being zero, see Fig.~\ref{fig:delta}(b).  Without the nonnegativity constraint, the content of the null space, in the case of LW, or excluded space, in the case of SVD, stays at zero.

One puzzling outcome of the tests so far is that the LW method and the SVD method, when its excluded space is identical to the null space, perform as well as~RL in restoring the 4-pixel ramp in Fig.~\ref{fig:4PixelRamp}, even though they struggle with the $\delta$-distribution.  The mystery is solved when looking at the null-space singular vectors in Fig.~\ref{fig:singulars} and the original ramp in Fig.~\ref{fig:4PixelRamp}: a ramp with an even number of pixels at the same intensity lacks any null-space contribution for our blurring function of~Eq.~\eqref{eq:t5}.  When the ramp's extension is changed to an odd number of pixels or the ramp is sloped, the LW and SVD methods can develop negative intensity values for the restored image, just like in the case of the $\delta$ distribution.  

\begin{figure}[!htb]
    \centering
    \begin{subfigure}[b]{0.46\linewidth}        
        \centering
        \includegraphics[width=\linewidth]{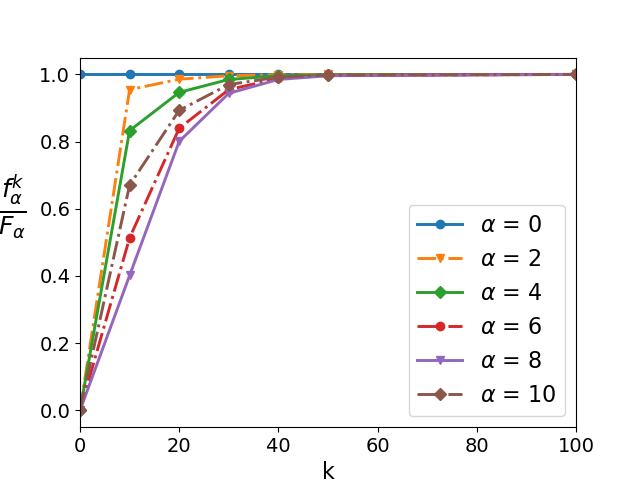}
        \caption*{(a)}
        \label{lambda1}
    \end{subfigure}
     \hspace*{1em}
    \begin{subfigure}[b]{0.46\linewidth}        
        \centering
    \includegraphics[width=\linewidth]{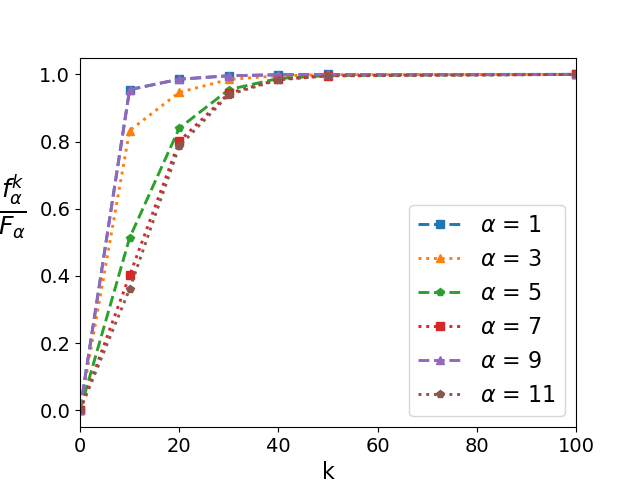}
        \caption*{(b)}
        \label{lambda2}
    \end{subfigure}
     \caption{Restored-image expansion coefficients in singular vectors, as a function of the RL iteration step $k$, normalized to the coefficients for the original $\delta$-distribution image, for even-index (a) and odd-index (b) singular vectors.}
    \label{fig:behavior}
\end{figure}

Given the benefits of the built-in nonnegativity constraint, we will further exclusively rely on the RL deblurring method.  Other methods may perform similarly when supplemented with this type of constraint.  In employing the RL method, the iteration progress is sometimes used as a tool.  Figure \ref{fig:behavior} shows how the expansion coefficients in singular vectors for the restored image behave as a function of the RL iteration step.  Here, the original image is the $\delta$-distribution, and the starting image is the uniform distribution.  One can observe that the higher the index of the singular vector, the generally slower the approach of the coefficient to its asymptotic value.  This can be understood in terms of the weakening impact of expansion coefficients on the amplification factors \eqref{eq:RLamplification} when the singular value decreases, as
\begin{equation}
    \frac{\text{d}A_i^{(k)}}{\text{d}f_\alpha^{(k)}} \propto \sigma_\alpha \, ,
\end{equation}
for $g_i >0$.  Null-space and low-singular-value coefficients will be affected by the positivity constraint and/or regularization, and in Fig.~\ref{fig:behavior} we observe a lumping in the pace of their approach to the asymptotic values.

\section{Exploring the limits of deblurring}
In the context of deblurring applications in quantitative research in nuclear and particle physics, it is important to understand when and how deblurring can falsify results.  We have observed the role that positive definiteness could play in restoring the null-space contributions to a distribution.  However, intuitively, one might expect the positive definiteness to play a minor role when the original distribution is uniformly far from zero on the distribution's general scale of variation.  Below, we will compare the performance of the classic RL method in those different situations.  When there is no positive definiteness to fall back on, regularization can improve the performance of a deblurring method.  However, the regularization must be a compromise with some adverse effects, too, as may be evident from Eq.~\eqref{func:min_mod}.

\subsection{High vs Low Contrast}

For high-contrast original images, such as in Fig.~\ref{fig:delta}(a) and \ref{fig:4PixelRamp}, we have observed that the deblurring with the nonnegativity constraint can restore the null-space components from the blurred images.  Notably, for the nonnegativity constraint to be impactful, the number of pixels at low intensity compared to the rest must be comparable to the dimension of the null space.  We now confront that situation to the situation with low-contrast original images, where the intensity for most pixels varies within a small range relative to the maximal intensity.

\begin{figure}[thb]
    \centering 
    \begin{subfigure}[b]{0.46\linewidth}        
        \centering
        \includegraphics[width=\linewidth]{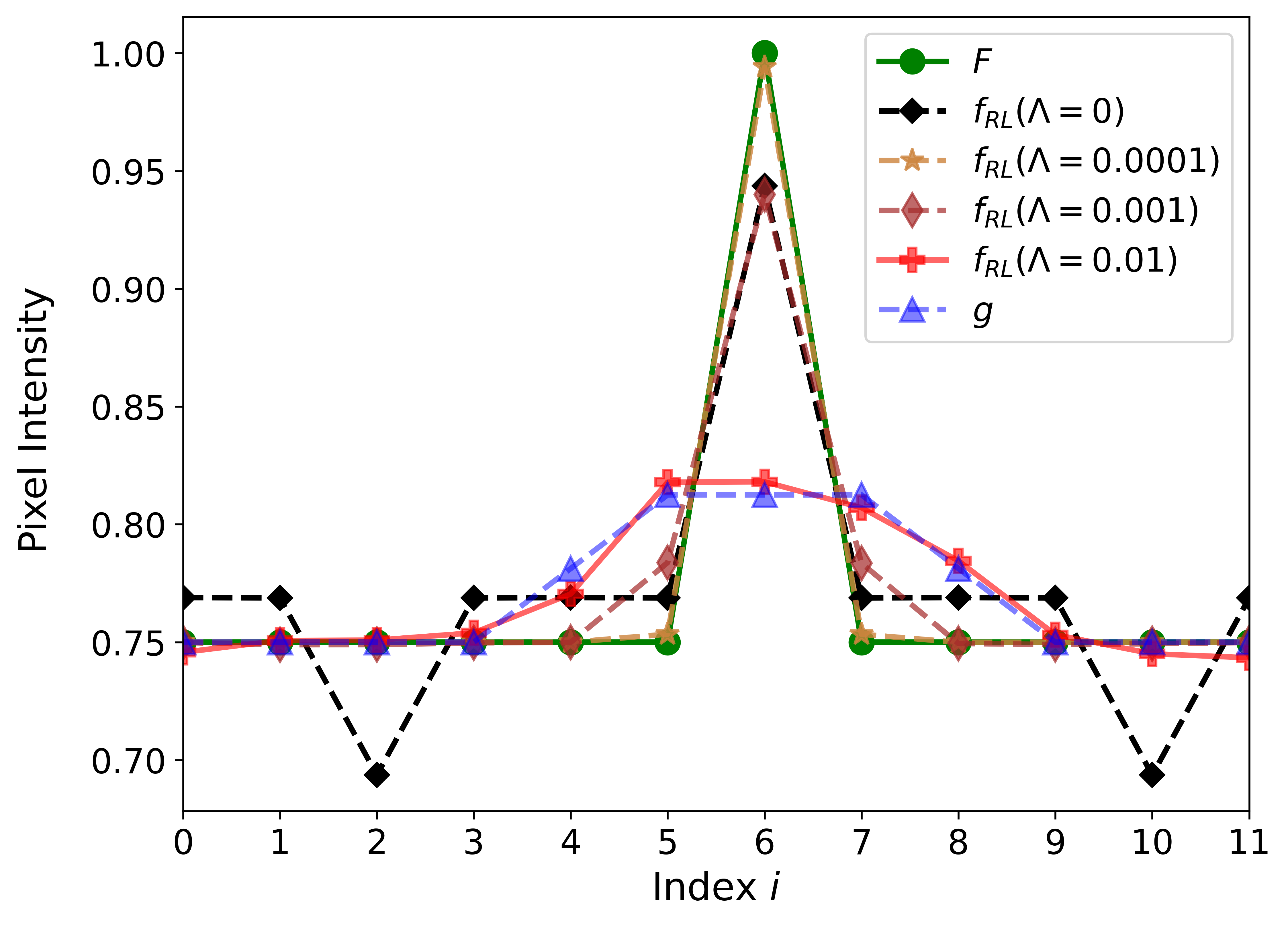}\\
        \caption*{(a)}
    \end{subfigure}
     \hspace*{1em}
    \begin{subfigure}[b]{0.46\linewidth}        
        \centering
    \includegraphics[width=\linewidth]{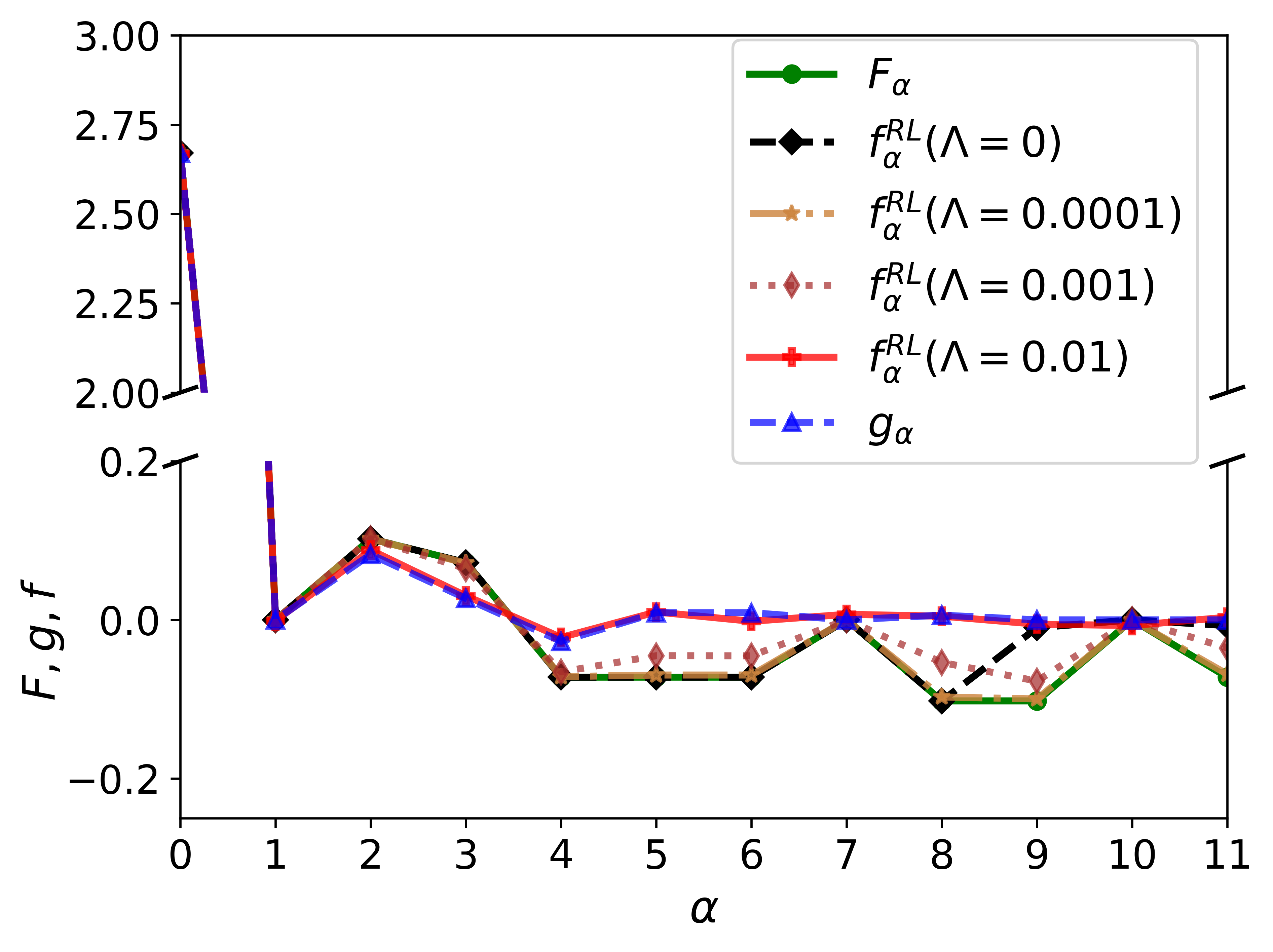}
        \caption*{(b)}
    \end{subfigure}
    \caption{Deblurring test with a single pixel with index 5 at intensity elevated by 0.25 over a high uniform background at 0.75.  \revised{(a)~Circles represent the original image $F$ and triangles - the image $g$ from blurring $F$ with the 5-pixel blurring function.  Moreover, the restored images~$f$ from different versions of the RL method are displayed.  Wide diamonds represent the image restored using the classical~RL method.  Stars, slim diamonds, and rectangles represent the images restored with the regularized RL method at $\Lambda = 0.0001$, 0.001, and 0.01, respectively.}  The~symbols corresponding to each case are connected with lines. (b)~Coefficients of decomposition in the basis of singular vectors for the cases in (a) vs singular-value index.}
    \label{fig:delta_LowContrast} 
 \end{figure}

An exemplary low-contrast case of a $\delta$-function superimposed on a high uniform background is illustrated in Fig.~\ref{fig:delta_LowContrast}(a).  When the classical RL method is applied to the blurred image, the method fails, on the scale of the $\delta$-function norm, in a similar manner as the methods without the nonnegativity constraint in the high-contrast case of Fig.~\ref{fig:delta}(a).  Again, we start the RL iterations in this illustration with a uniform image.  Although the nonnegativity is present in the RL method, it does not affect restoration when the background is sufficiently elevated.  Fig.~\ref{fig:delta_LowContrast}(b) complements the results of (a) with coefficients of decomposition in the basis of singular vectors, demonstrating that the coefficients for the null space are not restored for the classical RL method, unlike in Fig.~\ref{fig:delta}(b).  For impactful constraints, we turn to the regularized RL method, in which restoration is modified relative to the classical method based on differences in the intensity of adjacent pixels, thereby suppressing local extrema.  With a low $\Lambda = 0.0001$, the original image is satisfactorily restored in Fig.~\ref{fig:delta_LowContrast}(a).  \revised{The Root Mean Squared Error (RMSE) values for the restorations of the image in Fig.~\ref{fig:delta_LowContrast} and others, using different values of $\Lambda$, are provided in Table \ref{tab:RMSE}.}

\begin{table}[!thb]
\centering
\begin{tabular}{l l c}
\hline
\textbf{\textcolor{black}{Figure}} & \textbf{\textcolor{black}{Case}} & \textbf{\textcolor{black}{RMSE}} \\ 
\hline
\textcolor{black}{Fig.7} & \textcolor{black}{$f_{RL}(\Lambda=0)$} & \textcolor{black}{0.032558} \\
 & \textcolor{black}{$f_{RL}(\Lambda=0.0001)$} & \textcolor{black}{\textbf{0.002206}} \\
 & \textcolor{black}{$f_{RL}(\Lambda=0.001)$} & \textcolor{black}{0.022123} \\
 & \textcolor{black}{$f_{RL}(\Lambda=0.01)$} & \textcolor{black}{0.059650} \\
 & \textcolor{black}{$g$} & \textcolor{black}{0.061184} \\
\hline
\textcolor{black}{Fig.8(a)} & \textcolor{black}{$f_{RL}(\Lambda=0)$} & \textcolor{black}{0.033004} \\
 & \textcolor{black}{$f_{RL}(\Lambda=0.0001)$} & \textcolor{black}{\textbf{0.000377}} \\
 & \textcolor{black}{$f_{RL}(\Lambda=0.001)$} & \textcolor{black}{0.003743} \\
 & \textcolor{black}{$f_{RL}(\Lambda=0.01)$} & \textcolor{black}{0.037838} \\
 & \textcolor{black}{$g$} & \textcolor{black}{0.058463} \\
\hline
\textcolor{black}{Fig.8(c)} & \textcolor{black}{$f_{RL}(\Lambda=0)$} & \textcolor{black}{0.035454} \\
 & \textcolor{black}{$f_{RL}(\Lambda=0.0001)$} & \textcolor{black}{\textbf{0.000196}} \\
 & \textcolor{black}{$f_{RL}(\Lambda=0.001)$} & \textcolor{black}{0.002004} \\
 & \textcolor{black}{$f_{RL}(\Lambda=0.01)$} & \textcolor{black}{0.019756} \\
 & \textcolor{black}{$g$} & \textcolor{black}{0.057054} \\
\hline
\textcolor{black}{Fig.8(e)} & \textcolor{black}{$f_{RL}(\Lambda=0)$} & \textcolor{black}{0.066372} \\
 & \textcolor{black}{$f_{RL}(\Lambda=0.0001)$} & \textcolor{black}{\textbf{0.000171}} \\
 & \textcolor{black}{$f_{RL}(\Lambda=0.001)$} & \textcolor{black}{0.001903} \\
 & \textcolor{black}{$f_{RL}(\Lambda=0.01)$} & \textcolor{black}{0.017943} \\
 & \textcolor{black}{$g$} & \textcolor{black}{0.114109} \\
\hline
\end{tabular}
\caption{\textcolor{black}{RMSE values for the images of Figs.\ \ref{fig:delta_LowContrast} and \ref{fig:LowContrastExamples}, blurred ($g$) and restored ($f_{RL}$) at different regularization strengths $\Lambda$.}}
\label{tab:RMSE}
\end{table}

\begin{figure}[!hptb]
    \centering   
    \begin{subfigure}{0.44\textwidth}
        \includegraphics[width=\linewidth]{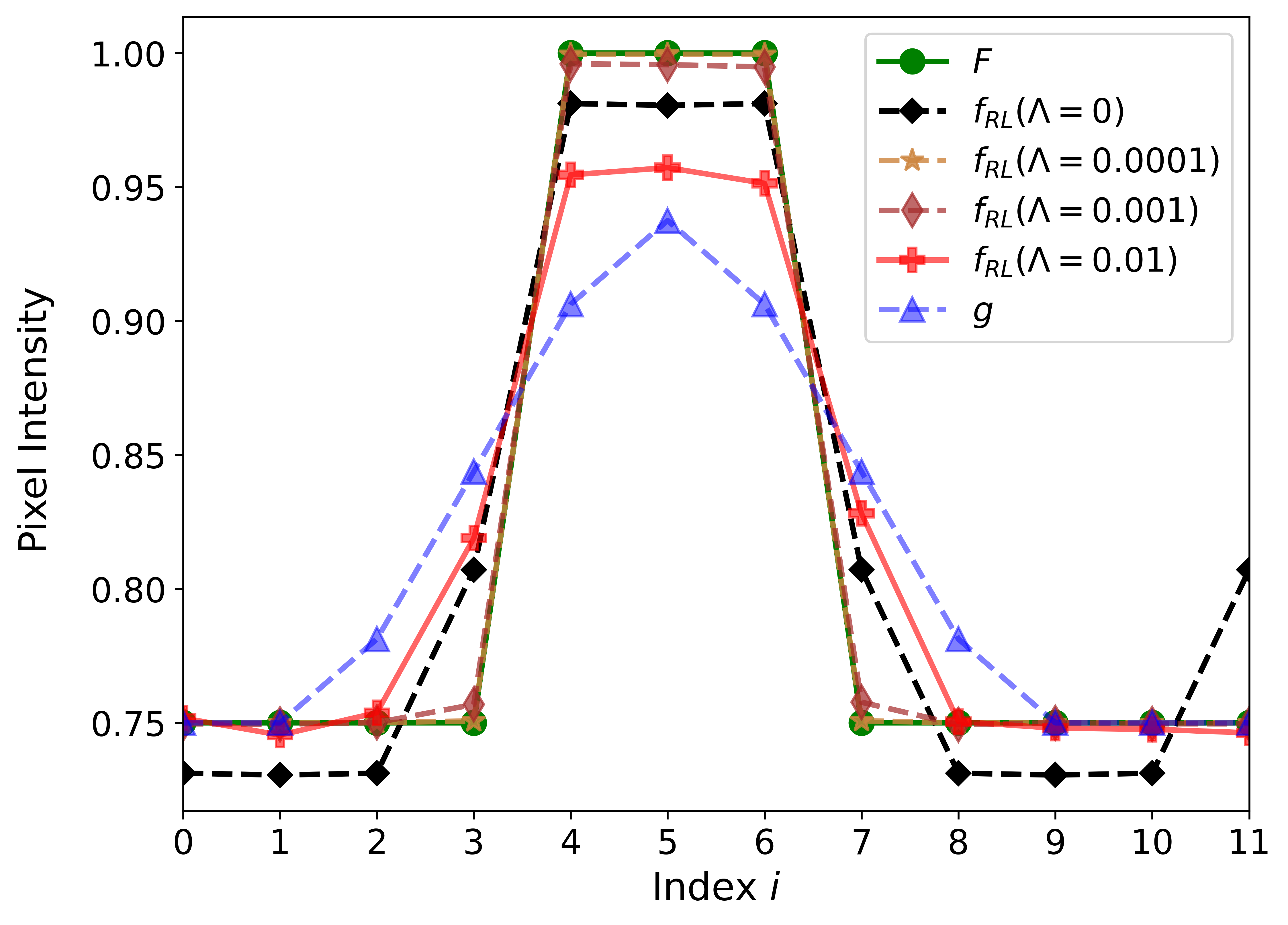}
         \caption*{(a)}
    \end{subfigure}
    \hspace*{1em}
    \begin{subfigure}{0.44\textwidth}
        \includegraphics[width=\linewidth]{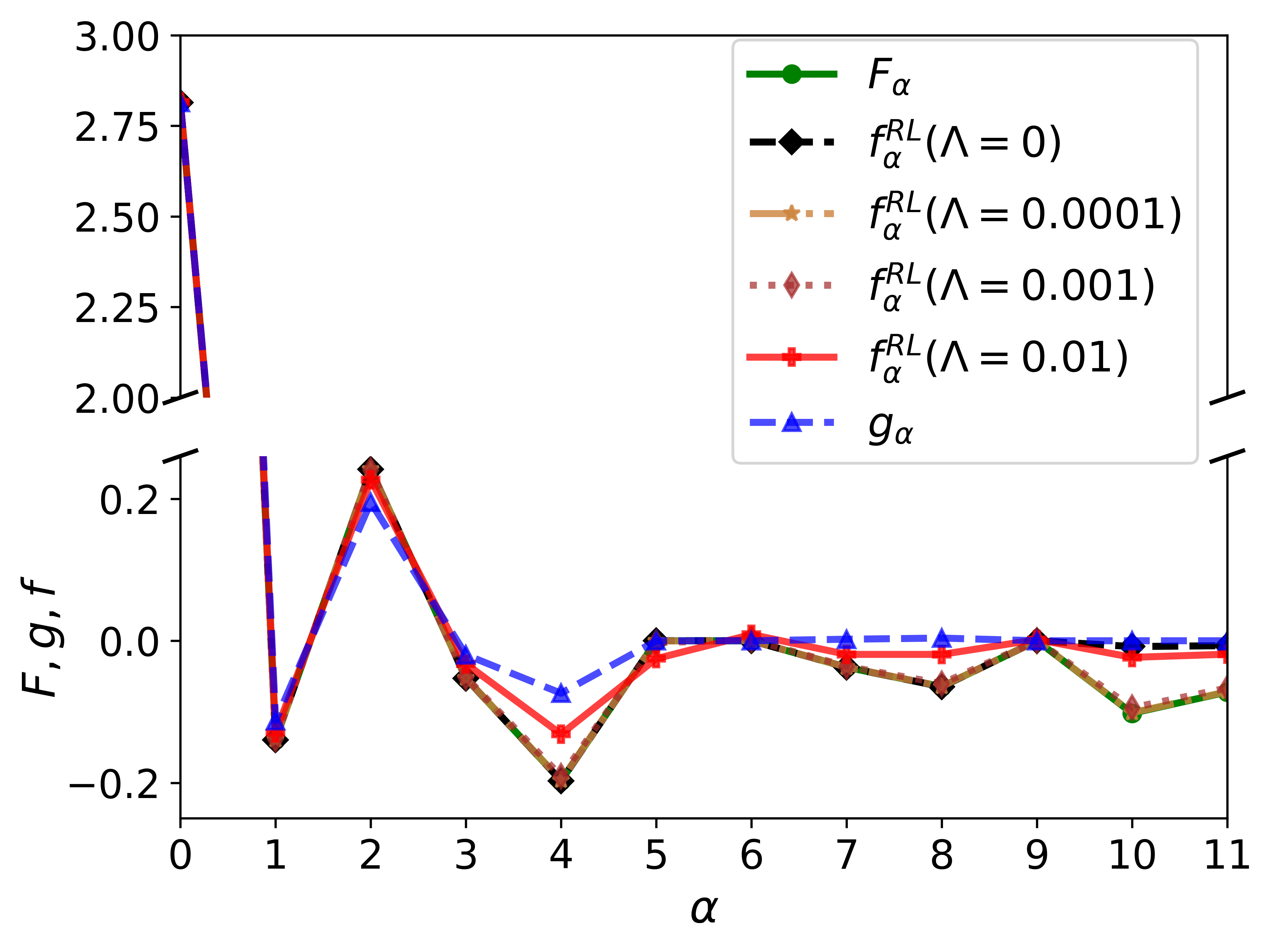}
         \caption*{(b)}
    \end{subfigure}  
        \vspace{0.4em} 
        \begin{subfigure}{0.44\textwidth}
        \includegraphics[width=\linewidth]{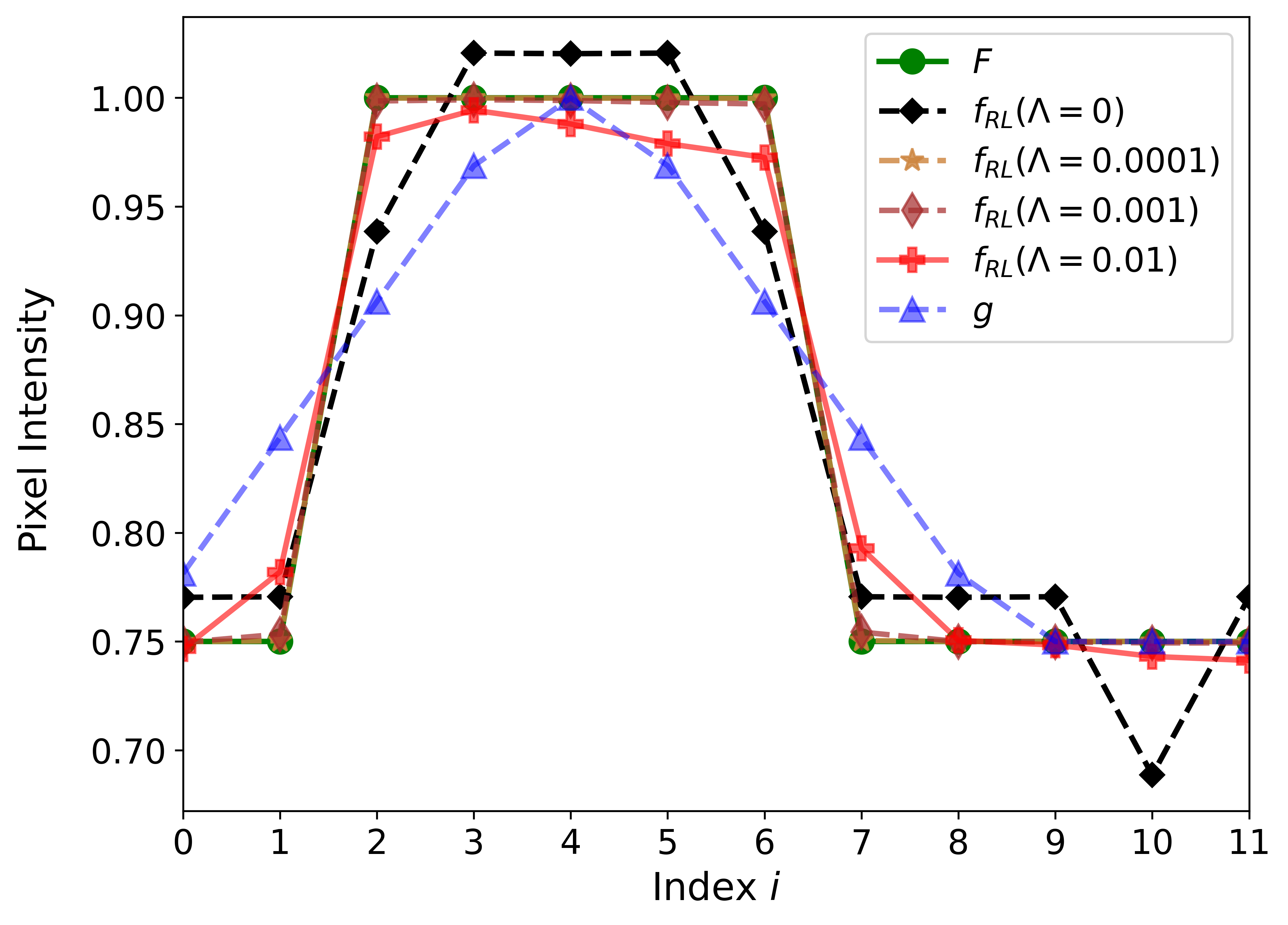}
         \caption*{(c)}
    \end{subfigure}
    \hspace*{1em}
    \begin{subfigure}{0.44\textwidth}
        \includegraphics[width=\linewidth]{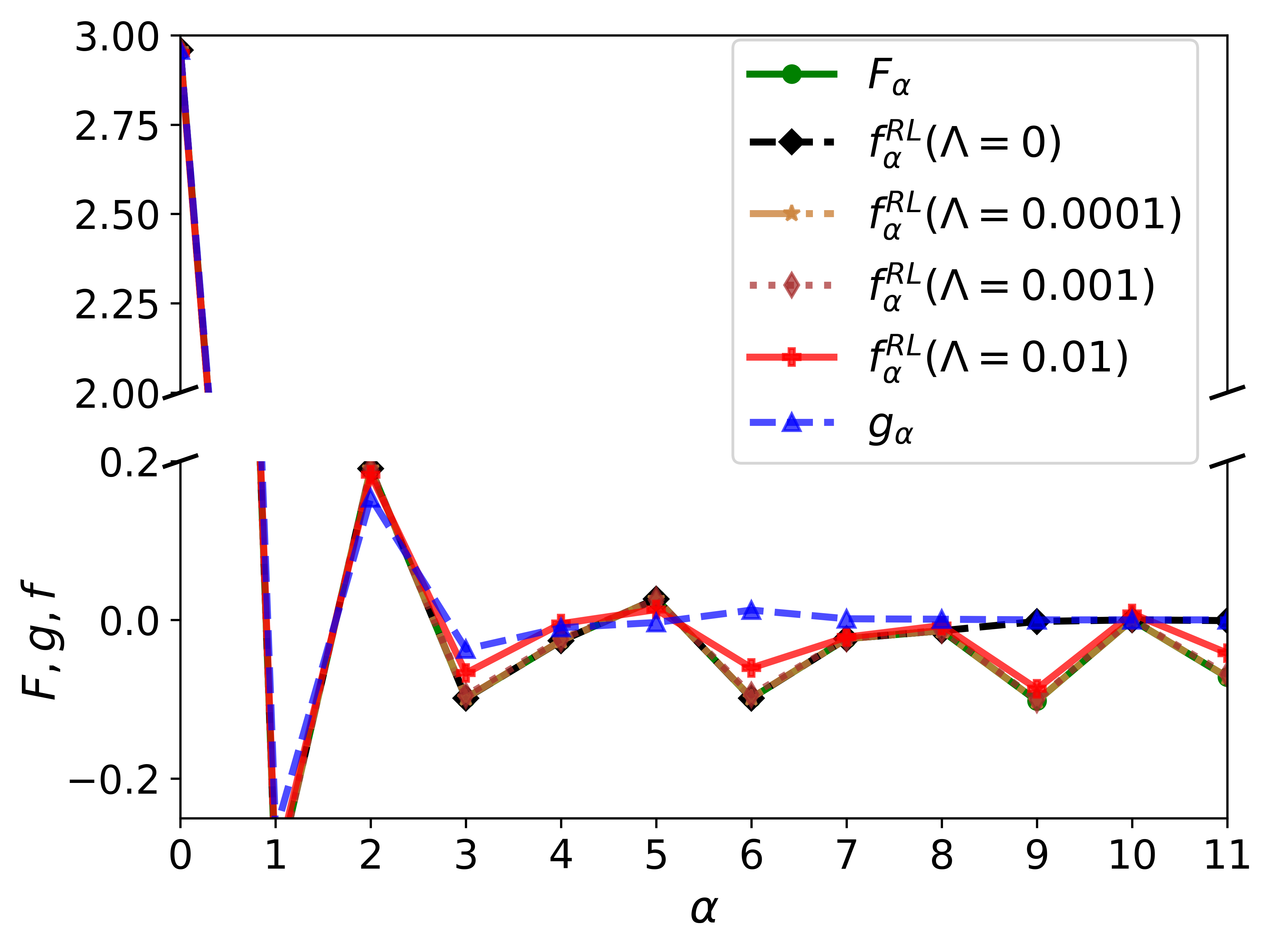}
        \caption*{(d)}
    \end{subfigure}
        \vspace{0.4em} 
        \begin{subfigure}{0.44\textwidth}
        \includegraphics[width=\linewidth]{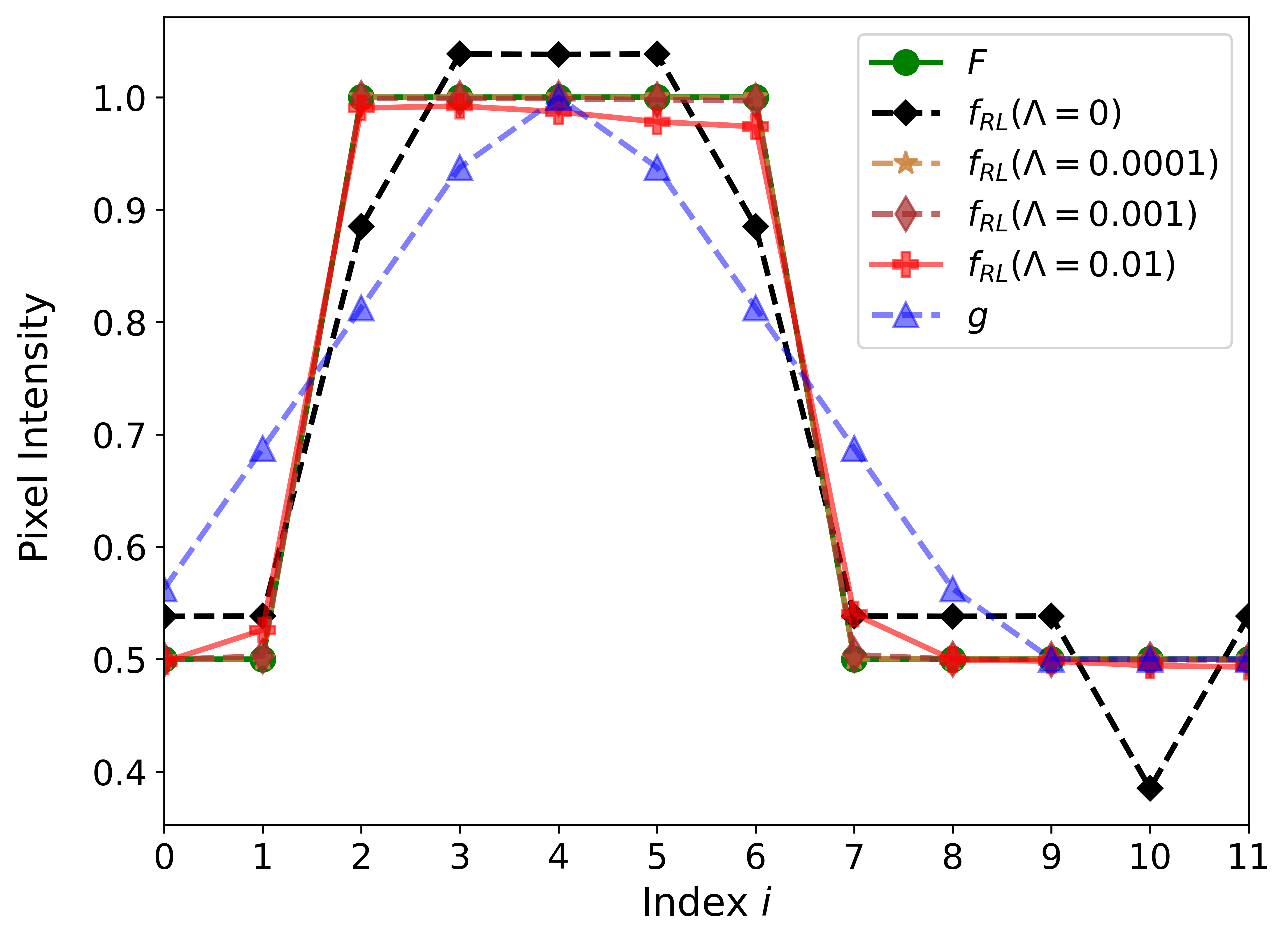}
         \caption*{(e)}
    \end{subfigure}
    \hspace*{1em}
    \begin{subfigure}{0.44\textwidth}
        \includegraphics[width=\linewidth]{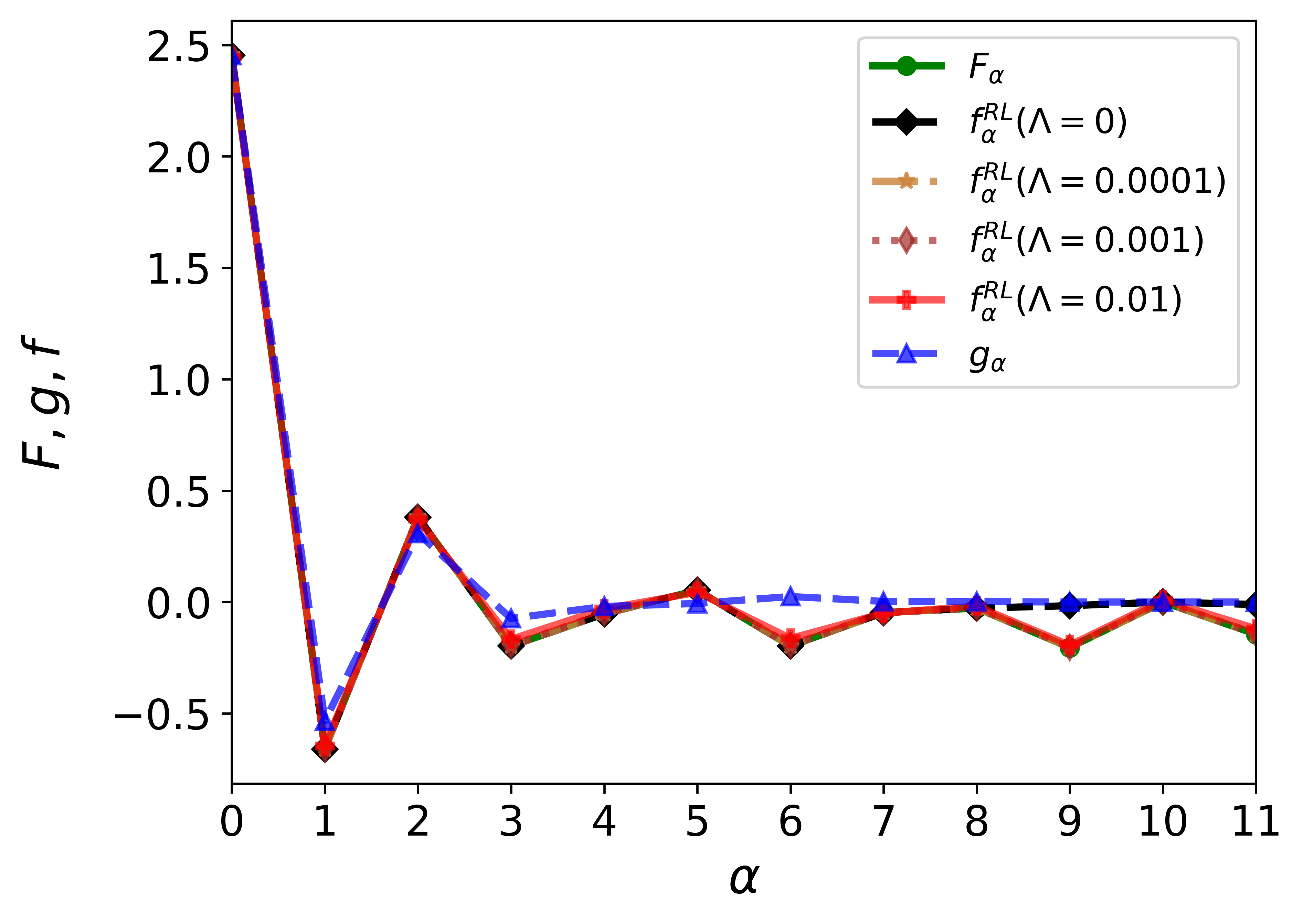}
        \caption*{(f)}
    \end{subfigure}
    \caption{Deblurring tests such as in Fig.~\ref{fig:delta_LowContrast} for other exemplary low-contrast images, with matching information for each of these presented in panels (a) and (b), (c) and (d), and (e) and (f), respectively.}           \label{fig:LowContrastExamples}
\end{figure}

Additional examples of image restoration at different contrast levels are shown in Fig.~\ref{fig:LowContrastExamples}.

\subsection{Impact of Initialization}

For a low-contrast original image $F$, it is possible to add a significant null-space component to that image, without violating intensity nonnegativity.  When such a combination is fed to the classical RL method as the starting guess $f^{(0)}$, mathematically, the~restoration iterations will never suppress the added spurious null-space component, returning that starting guess as the fully restored image $f$.  Even components added from low-singular-value vectors may only be slowly modified over iterations.

\begin{figure}[!htb]
    \centering
    \begin{subfigure}[b]{0.46\linewidth}  
        \centering
        \includegraphics[width=\linewidth]{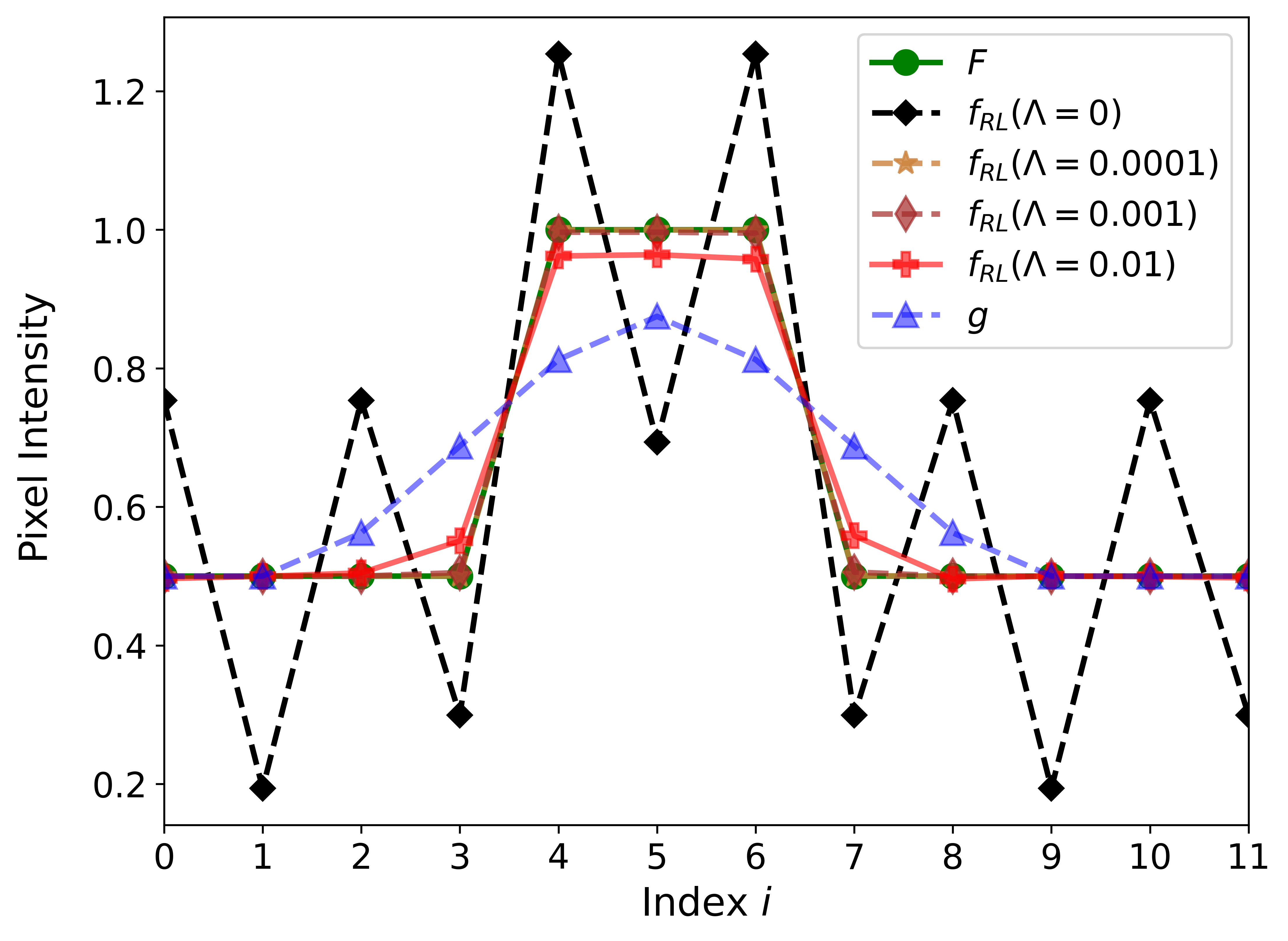}
        \caption*{(a)}
    \end{subfigure}
    \hspace*{1em}
    \begin{subfigure}[b]{0.46\linewidth}  
        \centering
        \includegraphics[width=\linewidth]{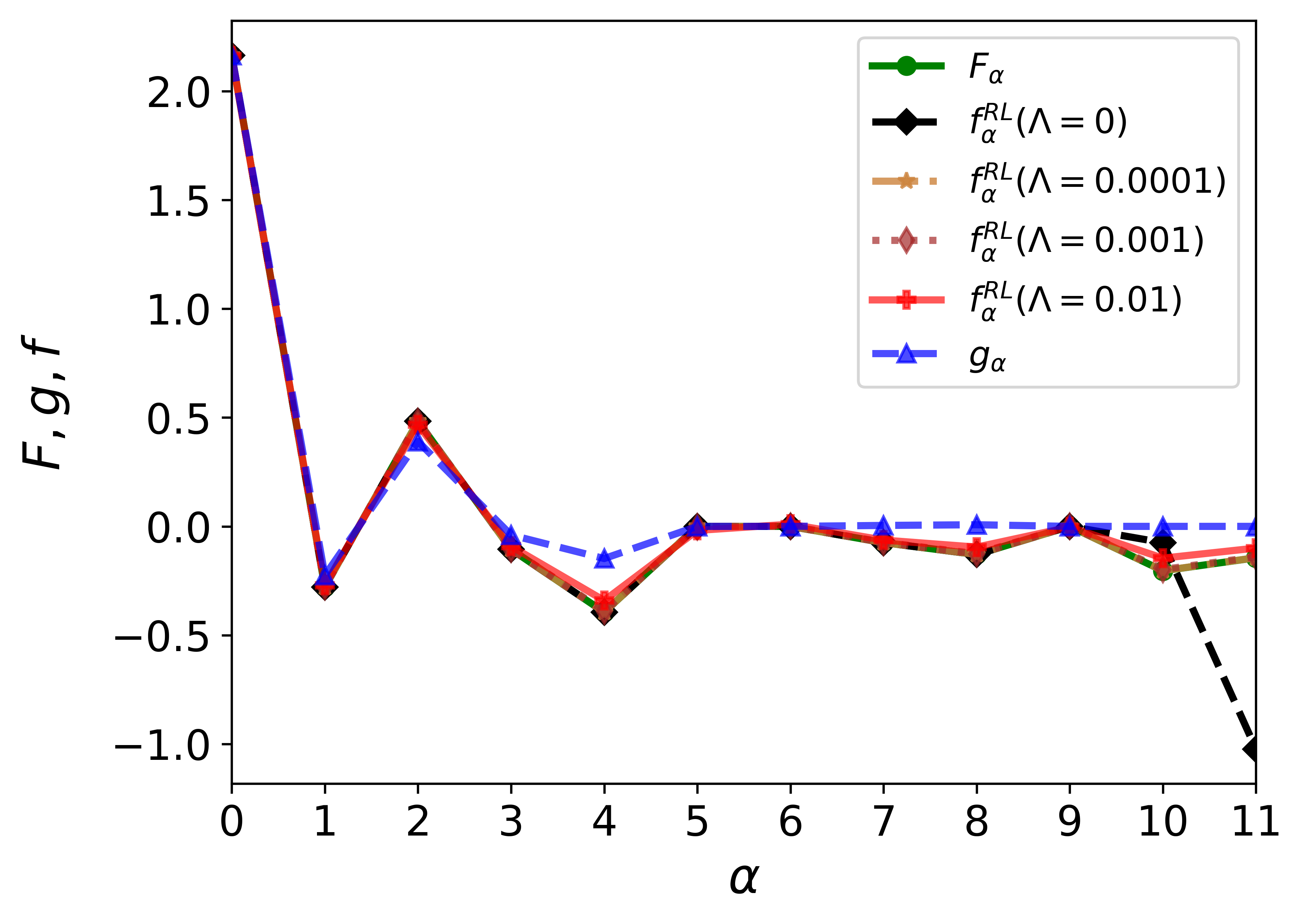}
        \caption*{(b)}
    \end{subfigure}
    
    \vskip 0.5cm  
    \begin{subfigure}[b]{0.46\linewidth}  
        \centering
        \includegraphics[width=\linewidth]{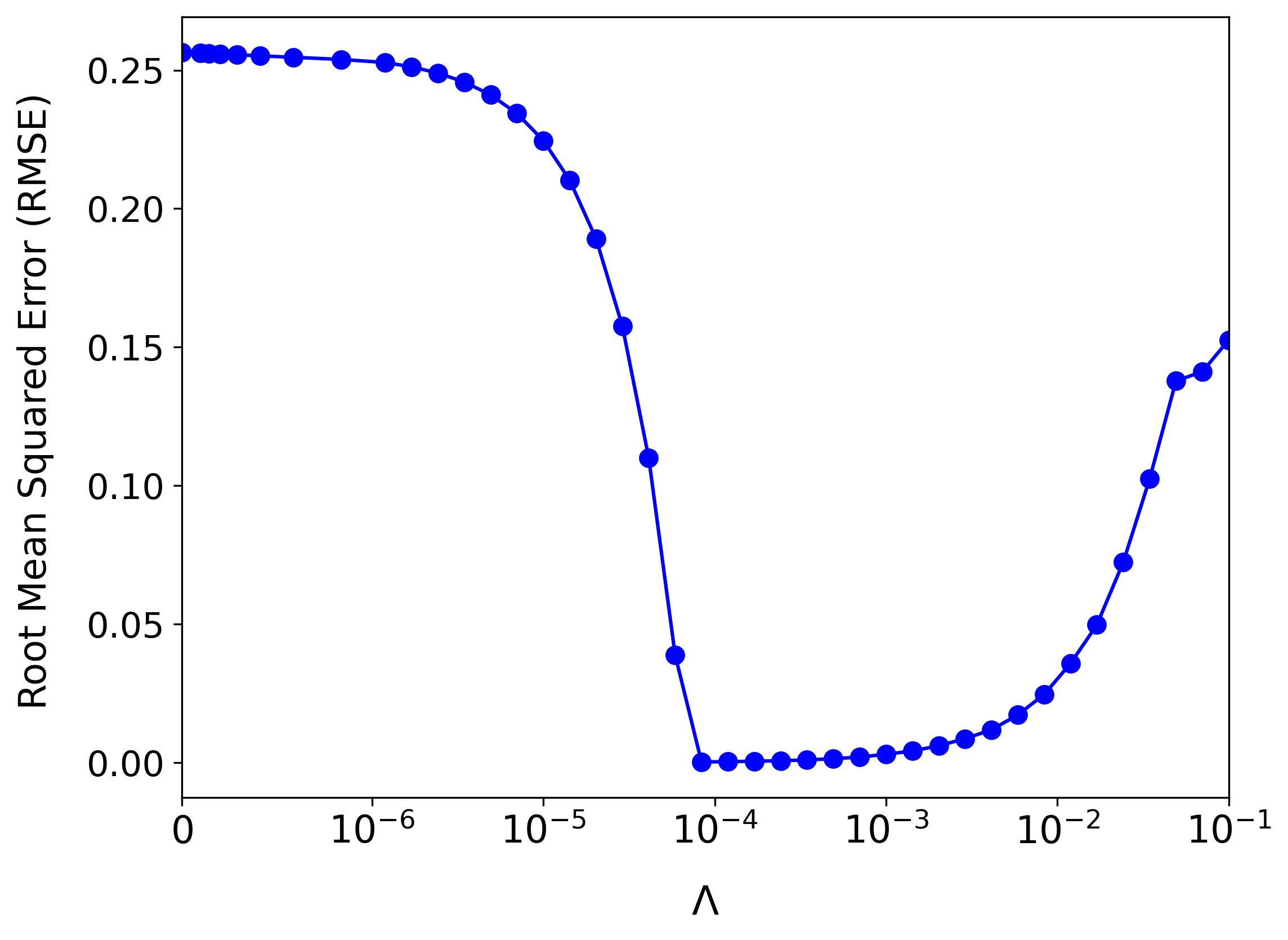}
        \caption*{(c)}
    \end{subfigure}
    
    \caption{Sensitivity of an RL restoration to null-space content at iteration start and to the regularization strength.  The RL iterations here start with the original image, with an $\alpha=11$ component added. Panels (a) and (b) show the images and singular-vector coefficients, respectively.   The classical RL method ($\Lambda = 0$) leaves the null-space content intact. An RL method with regularization ($\Lambda \gtrsim 10^{-4}$) suppresses the superimposed spurious null-space content, but excessive regularization, such as for $\Lambda=0.01$, can chip away the null-space content in the restored image compared to the original. Panel $(c)$ shows the RMSE values for the restored images as a function of $\Lambda$.   
    }
    \label{alpha_Odd3_mid}
\end{figure}

To avoid excessive spurious null-space or low-singular-value admixtures in the restored image, it helps to start the classical RL iterations with an image that lacks such components, such as a uniform image.  In Fig.~\ref{alpha_Odd3_mid}, we show what happens within the classical RL restoration when the restoration of a low-contrast image is started with a significant null-space component in the first guess for the image instead of the uniform guess.  The null-space component's strength persists at a significant level for $\Lambda = 0$ in the limit of many restoration steps.  Notably, the null-space components are unlikely to converge to the values of the original image in the absence of an impact from the nonnegativity constraint, regardless of how the iterations are started, including with the uniform image.  In Figs. In \ref{fig:delta_LowContrast} and \ref{fig:LowContrastExamples}, we see that classical restorations of low-contrast images that start with uniform guesses tend to underestimate the magnitudes of the null-space components in the restored images.  In particular,  the panels (a) and (b) in Fig.\ \ref{fig:LowContrastExamples} are analogous to Fig.~\ref{alpha_Odd3_mid} but differ importantly in initialization.  As shown in Figs.~\ref{fig:delta_LowContrast}, \ref{fig:LowContrastExamples}, and \ref{alpha_Odd3_mid}, the null-space components are brought under some control after a modest amount of regularization is applied.

\subsection{Impact of Regularization}

In examining the cases of restoration of interesting original images with few significant intensity jumps in Figs.~\ref{fig:delta_LowContrast} and \ref{fig:LowContrastExamples}, we can see that the regularization can dramatically improve image restoration.  The particular regularization of Eqs.~\eqref{tt} and \eqref{eq:I=} biases against multiple pronounced extrema and, with that, suppresses significant spurious null-space contributions that could be added to the restored low-contrast image for the classical RL method.

While a modest regularization within the RL method can help restore the null-space content of the original image, an excessive regularization can act to suppress not only that content but also contributions from singular vectors corresponding to low nonzero singular values.  This is illustrated in Figs.~\ref{fig:delta_LowContrast} and \ref{fig:LowContrastExamples}, which include results from RL restorations with progressively strengthening regularization.  As the regularization strength increases, the features of the restored images evolve from those of the original image toward those of the blurred image. Fig.~\ref{alpha_Odd3_mid}(c)  illustrates the dependence of RMSE for the restored image compared to the original image on the parameter $\Lambda$ for the regularization strength.  Excessively low values of $\Lambda$ may not be sufficient to tame the null-space components intentionally superimposed on the image starting guess.  The growth of RMSE with the growth of $\Lambda$ for the higher $\Lambda$ values can be seen, reflecting the gradual blurring in the features of the restored image.  In effect, there is an optimal range of $\Lambda$ in a restoration problem for low-contrast images.  The choice of regularization amounts to trading a stronger blurring in the image before deblurring for a weaker blurring afterward.

\subsection{\revised{More Pixels}}

\begin{figure}[!hptb]
    \centering   
    \begin{subfigure}{0.65\textwidth}
     \centering   
        \includegraphics[width=\linewidth]{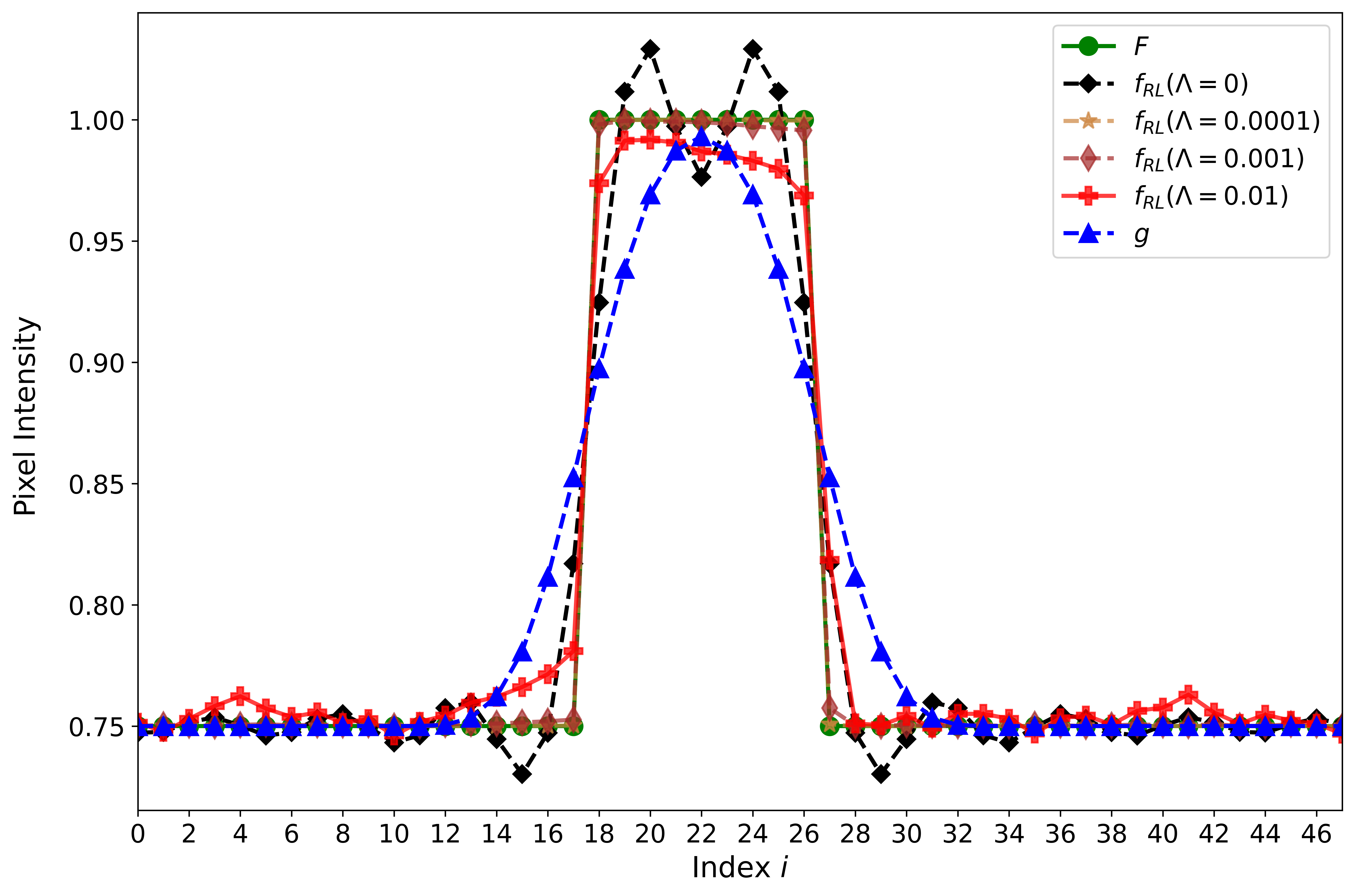}
         \caption*{(a)}
    \end{subfigure}
    \vspace*{1ex}
    
    \begin{subfigure}{0.65\textwidth}
     \centering   
        \includegraphics[width=\linewidth]{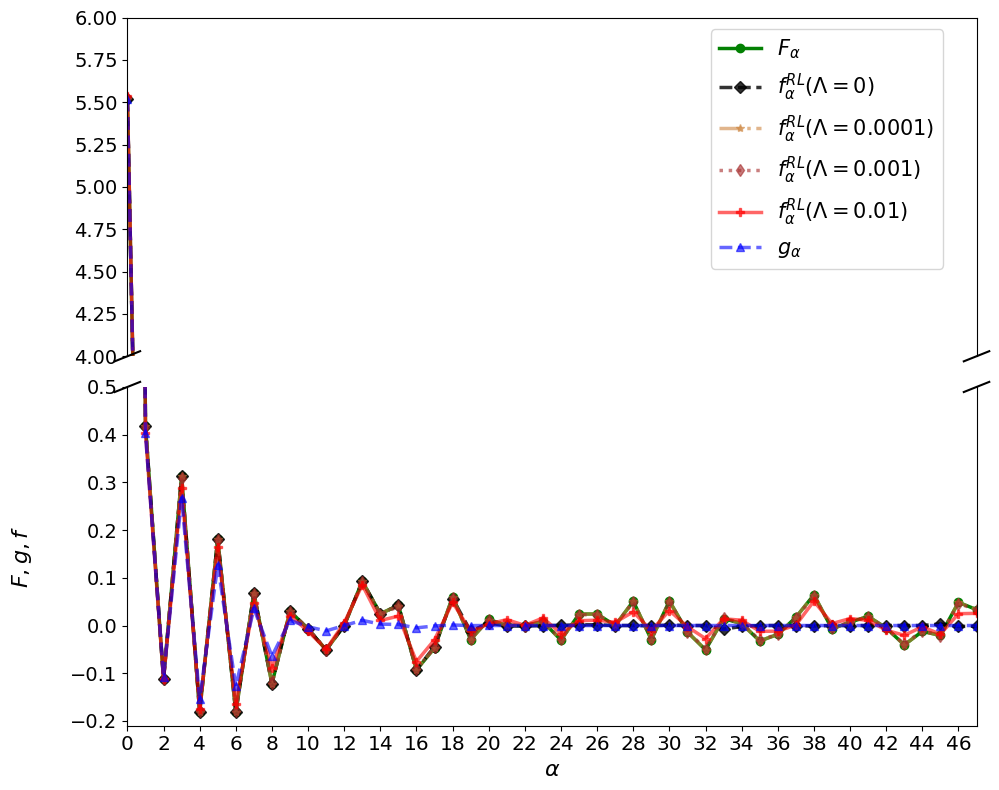}
         \caption*{(b)}
    \end{subfigure}  
         \caption{\revised{Deblurring tests for 48 pixels and the original low-contrast image similar to those represented in panels (a) and (b) of Fig.~\ref{fig:LowContrastExamples}.}  }         \label{fig:48Pixels}
\end{figure}

\revised{The case of just 12 pixels and a 5-pixel blurring function allows for an exploration of details in blurring and deblurring.  For more pixels and wider blurring functions, the features need to be looked at in a more statistical manner.  The separation between the null space and its complement may become less clear-cut, as the singular values may fall by many orders of magnitude with increasing index before hitting zero.  Even with the need for some perspective shift, our experimentation reveals similarities in the deblurring across both higher and lower pixel counts.}
\revised{
As an example, we take the case of 48 pixels and start with the 5-bin blurring function of Eq.~\eqref{eq:t5}.  We construct a wider 13-bin blurring function, commensurate with the larger pixel space, by taking a third power of the 5-bin matrix, which is equivalent to a trifold blurring with the 5-bin function.  With this, the singular vectors remain the same as for the 5-bin function, but now in the 48-pixel space, whereas the singular values are raised to the third power, causing the spectrum to fall off by orders of magnitude with the index before reaching the null space.  In Fig.~\ref{fig:48Pixels}, we show the case of blurring and deblurring of a low-contrast ramp similar to the one for 12 pixels in Fig.~\ref{fig:LowContrastExamples}(a).  Again, the restoration with no regularization yields a much distorted image. A relatively wide range of regularization strengths yields satisfactory restored images, and only at $\lambda \sim 0.01$ do they become excessively biased. }

\section{From Experimentation to Data Analysis}

We briefly discuss how this current backbone experimentation connects to the parallel nuclear data analysis. The blurring due to reaction-plane resolution in nuclear collisions \cite{danielewicz_deblurring_2022}, i.e., misassigning the azimuthal angle relative to the reaction plane to an emitted particle, nominally occurs in three dimensions.  However, this blurring and deblurring occur only in one dimension of the azimuthal angle around the collision beam axis.  The deblurring is repeated for separate velocity bins along and perpendicular to the beam.  In other recent nuclear applications in nuclear physics, such as decay-energy and $\gamma$-ray spectra~\cite{nzabahimana_deconvoluting_2023, xu_reconstruction_2024} and the imaging of spatial sources from particle correlations~\cite{tam_source_2025, nzabahimana_source_2023}, deblurring is only one-dimensional.  However, the applications may progress to more dimensions when, e.g., three-dimensional features of the sources are sought.

In the reaction-plane deblurring~\cite{danielewicz_deblurring_2022}, the particle yield distributions are periodic, like here, as the azimuthal angle is used in the process.  The transfer matrix for decay-energy spectrum deblurring is constructed by simulating the detection system's performance and is independent of the distribution it operates on, as here.  In the source-imaging from correlations, the transfer matrix is constructed theoretically and may be folded with detector resolution.   It is independent of the source that is inferred.  The case of reaction-plane deblurring is an interesting situation in which the intrinsic consistency of the deblurring procedures can be tested as they are applied in regions where detector efficiency nears 100\% and falls towards zero.

Elaborating on the reaction-plane case, the detectors perform well at detecting charged particles moving at high velocities, such as those close to the target velocity, but perform poorly at low velocities, approaching the target at rest in velocity space.  When only the forward-velocity region is analyzed, exploiting the central limit theorem, the transfer matrix can be well estimated from a few moments of the detected distribution~\cite{voloshin_flow_1996}, making it independent of the deblurring outcome.  With physically expected reflection symmetry relative to the reaction plane, the transfer matrix in the ideal-detector limit acquires a straightforward Fourier SVD form even ahead of any angular discretization:
\begin{equation}
\label{eq:Fourier}
\frac{d N}{d \phi'} = \int d\phi \, T(\phi'|\phi) \, \frac{d N}{d \phi} \simeq \sum_n \sigma_n \, \cos{n \phi'}  \int d\phi \,\cos{n \phi}  \, \frac{d N}{d \phi} \, .
\end{equation}
Here, $\phi$ and $\phi'$ are angles relative to the actual and estimated reaction plane, respectively, and $N$ is particle yield at a particular longitudinal velocity and transverse velocity magnitude.  The singular values are $\sigma_n = \langle \cos{n \Delta \Phi} \rangle$, where $\Delta \Phi$ is the deviation in the angle of the estimated plane from the actual.  However, the low-velocity transfer matrix results from a complex interplay of detector performance, especially when azimuthally asymmetric, and the inferred distribution.  Still, when the projectile is the same or nearly the same as the target, the distributions from deblurring in the velocity space near the projectile and target should physically be near mirror reflections of each other.  This can serve as the mentioned test of deblurring validity, as the procedures' details and complexity differ across regions, with the singular values in the SVD decomposition considered here never approaching unity in the low-velocity areas.  Notably, before deblurring, the collision data have been quantified in terms of Fourier coefficients, as evident in \eqref{eq:Fourier}, primarily $n=1$ and 2.

Deep Learning might be used to image nuclear collisions, but there are currently obstacles in generating training data for the networks, no matter how the training is envisioned.  The theoretical transport models of collisions are phenomenological and may miss essential features of the original images one seeks to access. They are also computationally expensive, making it difficult to reach the statistics needed for training. The situation may change once the essential features of the images have been settled, allowing for simplifying the training.

\section{Summary}

We examined the iterative RL deblurring method and partially other methods in the context of the utility of deblurring in nuclear and high-energy physics.  We ignore explicit noise in the images for the time.  SVD of the blurring matrix is an essential tool in understanding both the blurring and deblurring processes.  The singular vectors generally exhibit a faster variation with pixel position for decreasing singular values.  Symmetries of the blurring matrix generally give rise to symmetries of the singular vectors and may make the left- and right-hand vectors coincide. Blurring strength may be characterized by the falloff of singular values with a singular-vector index.  A sufficiently strong blurring can give rise to a null space in the SVD.  A one-dimensional model system of 12 pixels with periodic boundary conditions has served to illustrate our points.  For high-contrast images, where the number of pixels at zero or near-zero intensity relative to the image's maximal intensity is comparable to the null-space dimension, the nonnegativity constraint built into a deblurring method can help restore the null-space content in the processed image.  For lower-contrast images, a regularization term that penalizes multiple extrema in the restored image helps control the null-space content.  However, overly strong regularization may start affecting the restored image, such as blurring, which one is trying to correct. Despite the simplicity of the primary model we studied here, it can come close to the data analyses carried out in parallel in various respects. 

\revised{For a clearer picture, a similar analysis of deblurring performance in the presence of noise in the blurred images is underway.}

\printbibliography

\end{document}